\documentclass[10pt,a4paper]{article}
\usepackage[utf8]{inputenc}
\usepackage[english]{babel}
\usepackage{amsmath}
\usepackage{amsfonts}
\usepackage{amssymb}
\usepackage{tabularx,multirow}
\usepackage{graphicx}
\usepackage{subfigure}
\usepackage{mathrsfs}
\usepackage{mathenv}
\usepackage{dsfont}
\usepackage{mathtools}
\usepackage[usenames, dvipsnames]{xcolor}
\usepackage{pgf,tikz}
\usetikzlibrary{arrows}
\usepackage[left=2cm,right=2cm,top=2cm,bottom=2cm]{geometry}

\title{Variance reduction method for particle transport equation in spherical geometry}
\author{X. Blanc$^1$, C. Bordin$^2$, G. Kluth$^2$, G. Samba$^2$,\\ \\
\footnotesize $^1$ Univ. Paris Diderot, Sorbonne Paris Cité, \\
\footnotesize Laboratoire Jacques-Louis Lions, UMR 7598, UPMC, CNRS,\\
\footnotesize F-75205 Paris, France \\ \\
\footnotesize$^2$ CEA, DAM, DIF, 91297 Arpajon Cedex, FRANCE
}
\date{\today}
\def\bx{\mathbf{x}}
\def\bOmega{\mathbf{\Omega}}
\def\RR{\mathbb{R}}
\def\Ei{\operatorname{Ei}}
\def\red{\textcolor{red}}

\newtheorem{remarque}{Remark}[section]

\numberwithin{equation}{section} 
\usetikzlibrary{arrows}
\begin{document}
\maketitle

\begin{abstract}
  This article is devoted to the design of importance sampling method for the Monte Carlo simulation of a linear
  transport equation. This model is of great importance in the simulation of inertial confinement fusion
  experiments. Our method is restricted to a spherically symmetric idealized design : an outer sphere emitting 
  radiation towards an inner sphere, which in practice should be thought of as the hohlraum and the fusion
  capsule, respectively. We compute the importance function as the solution of the corresponding stationary adjoint problem. Doing
  so, we have an important reduction of the variance (by a factor 50 to 100), with a moderate increase of
  computational cost (by a factor 2 to 8). 
\end{abstract}

\section{Introduction}
\begin{figure}
  \centering
  \includegraphics[width=0.7\textwidth]{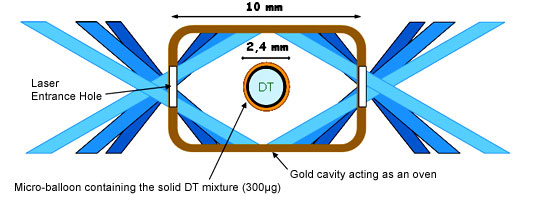}
  \caption{Schematic view of the Hohlraum and the target}
  \label{fig:1}
\end{figure}

\begin{figure}
  \centering
  \includegraphics[width=\textwidth]{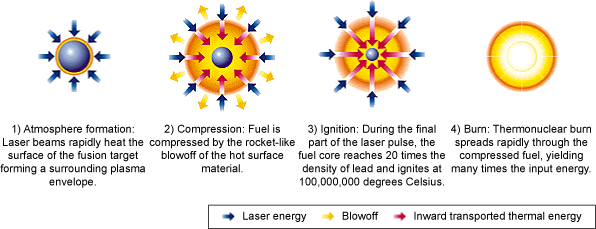}
  \caption{The concept of ICF (inertial confinement fusion) taken from http://www.lanl.gov/projects/dense-plasma-theory/background/dense-laboratory-plasmas.php}
  \label{fig:2}
\end{figure}

In inertial confinement fusion (ICF) experiments, a small ball of hydrogen (the target) is submitted to intense
radiation by laser beams. These laser beams are either pointed directly to the target (direct drive approach), or pointed
to gold walls of a hohlraum in which the target is located (indirect drive approach, see Figure~\ref{fig:1}). These gold walls heat up, emitting X-rays toward the
target. The outer layers of the target are heated up, hence ablated. By momentum conservation, the inner part of
the target implodes (this is usually called the rocket effect). Hence, the pressure and temperature of the hydrogen
inside the target increase, hopefully reaching the thermodynamical conditions for nuclear fusion. This process is
summarized in Figure~\ref{fig:2}.

The numerical simulation of such an experiment involves many physical phenomena such as hydrodynamics, radiation
transfer, neutronics, etc... In the present article, we focus on the simulation of radiation, that is, the
transmission of the (X-ray) energy to the target. A simplified model for this is the grey radiative transfer
equation:
\begin{equation}
\left\{
\begin{aligned}
\partial_t u + \bOmega\cdot \nabla u + \kappa_t u  & = \kappa_s \int_{S^2}  u(t,\bx,\bOmega') 
k(\bx,\bOmega',\bOmega)d\boldsymbol{\Omega'} + Q(\bx) \\
	 u(t=0,\bx,\bOmega) &= g(\bx,\bOmega),
\end{aligned}
\right.
\label{E_transport}
\end{equation}
where the solution $u$ is the radiation intensity and depends on the time $t$, the position $\bx\in
\mathbb{R}^3$, the direction of propagation $\bOmega \in S^2 $. The term $Q(\bx)$
represents a source of radiation. In the present case, $Q(\bx)$ is a modelling of the emission of X-rays by the
hohlraum walls. Furthermore, $\kappa_t$ is the total 
cross-section. It satisfies $\kappa_t = \kappa_a +
\kappa_s$, where $\kappa_a\geq 0$ is the absorption
cross-section and $\kappa_s\geq 0$ the scattering cross-section. The kernel $k(\bx,\bOmega',\bOmega)$ is a
probability density with respect to $\bOmega'$ and $\bOmega$, that is, $k\geq 0$ and $\displaystyle \int
k(\bx,\bOmega',\bOmega)d\bOmega' = \int
k(\bx,\bOmega',\bOmega)d\bOmega=1$. Note that we have assumed here that we use units such that the speed of light
is $c=1$.

Equation (\ref{E_transport}) may be simulated using a Monte Carlo method. If so, the probability distribution $k$
may be interpreted as the probability density associated to the new direction propagation $\bOmega$
for a particle having a shock with initial direction $\bOmega'$. 
In Monte Carlo simulations of such situations, variance reduction
methods are important to reduce the statistical noise. Indeed, as the target implodes, hydrodynamic instabilities
develop, which are a source of energy loss. Should this loss be too important, the experiment would be
compromised. Thus, it is important to have a precise numerical description of these instabilities. In the case of
Monte Carlo simulations, this implies a statistical noise as small as possible (at least smaller than the amplitude of
the instabilites). A small variance is particularly important on the target boundary. 

A widely used reduction variance technique in such a situation is the importance sampling method. \red{It may be
  summarized as follows:
  \begin{enumerate}
  \item\label{item:1} Calculate the importance function (in our case, the solution of the adjoint equation);
  \item\label{item:2} Use the importance function to modify the transport equation, and apply a Monte Carlo method;
  \item Calculate the forward intensity from \ref{item:1} and \ref{item:2}
  \end{enumerate}
} 
\red{Importance sampling is a well-known reduction variance method, which has
  been applied to transport problems in many situations. We refer for instance to the textbooks \cite{lux-koblinger}
  \cite{spanier-gelbard} for a general presentation. The key-point in such a method is the way one computes the
  importance function. If it is solution to the adjoint problem, then one achieves a zero-variance method. However,
solving the adjoint problem is at least as difficult as solving the direct problem at hand. Therefore, many methods
using approximations of the adjoint solution have been developped. This is the spirit of the exponential transform
(see \cite{depinay-these} and \cite{lux-koblinger}). In some situations, a diffusion approximation is used for this
calculation, as for instance in \cite{wright-etal}. In other situations, discrete ordinates approximation is
preferred \cite{wagner-haghighat}. The method which is the closest to the one presented here is probably
\cite{bal-davis-langmore}, in which the adjoint equation is formulated as an integral equation, and solved using a
space discretization. An importance difference is, however, that when solving the adjoint problem, the scattering
is neglected in \cite{bal-davis-langmore}. Here, we use the same kind of method, but taking advantage of the
radially symmetric geometry, we are able to take scattering effects into account.}

\medskip

The article is organized as follows: in Section~\ref{sec:monte-carlo-method}, we give a rapid presentation of the
Monte Carlo method applied to transport equations, then of the importance sampling method. This method is based on
the computation of an importance function, which is the subject of Section~\ref{sec:comp-import-funct}. In
Section~\ref{sec:numerical-results}, we present some numerical experiments, while the appendices contain some
technical result which we do not want to detail in the main body of the article.

\section{Monte Carlo method for transport equations}
\label{sec:monte-carlo-method}

\subsection{Natural method}
\label{sec:natural-method}

We give in this subsection a short overview of the application of Monte Carlo method applied to transport
equation. More details and mathematical justifications are given in \cite{lapeyre-pardoux-sentis}. We
concentrate here on practical aspects. 

Considering equation \eqref{E_transport}, we define the following quantities (we assume here that $Q$ and $g$
are integrable functions):
\begin{equation}\label{eq:alpha}
  \alpha = \int_{S^2}\int_{\mathcal D} g(\bx,\bOmega)d\bx d\bOmega, \quad \overline g(\bx,\bOmega) = \frac 1 \alpha
  g(\bx,\bOmega), 
\end{equation}
\begin{equation}\label{eq:beta}
  \beta = \int_{S^2}\int_{\mathcal D} Q(\bx)d\bx d\bOmega = 4\pi \int_{\mathcal D} Q(\bx)d\bx,
  \quad \overline Q(\bx) = \frac 1 \beta Q(\bx)
\end{equation}
Here, ${\mathcal D}$ is the spatial domain. Hence, $\overline g$ and $\overline Q$ are probability measures on the
phase space ${\mathcal D}\times S^2$. Of course, we assume that both $g$ and $Q$ are non-negative, which is
physically relevant. 

We first deal with the case $Q = 0$, and then extend it to the general case. 
We define $N$ independent realizations $(X_i(t),\Omega_i(t))$ of the jump Markov process $(X(t),\Omega(t))$ as
follows:
\begin{enumerate}
\item $(X_i(0),\Omega_i(0))$ are drawn independently of each other, following the law $\overline g(\bx,\bOmega)d\bx
  d\bOmega$. \\ To each of
  them is assigned a weight $w_i(0)=\frac 1 N$. 
\item Between jumps, $(X_i,\Omega_i)$ follows the characteristics of Equation \eqref{E_transport}, that is,
  \begin{displaymath}
    \left\{
      \begin{aligned}
        &\dot X_i(t) = \Omega_i(t), \\
        &\dot \Omega_i(t) = 0,
      \end{aligned}
\right.
  \end{displaymath}
which is equivalent to the fact that $\Omega_i$ is constant\red{\footnote{\red{Note that in curvilinear coordinates, this is not the
  case. For instance, in Section~\ref{sec:spherical-case} below, spherical coordinates are used, hence the
  direction $\mu_i$ is a non-trivial function of $t$.}}} and 
\begin{equation}\label{eq:trajecto}
  X_i(t) = X_i(t_0) + (t-t_0)\Omega_i. 
\end{equation}
Moreover, the weight $w_i(t)$ is assumed to satisfy the equation $\dot w_i(t) + (\kappa_t - \kappa_s)w_i(t) = 0$,
that is,
\begin{equation}\label{eq:amortissement}
  w_i(t) = w_i(t_0) e^{-\left(\kappa_t - \kappa_s\right)(t-t_0)}.
\end{equation}
\item Time jumps are defined by a Poisson process ${\mathcal N}(t)$ of intensity one as follows: if the process
  ${\mathcal N}(t\kappa_s)$ has a jump at time $t$, then $\Omega_i(t)$ has a jump, and the conditional law of
  $\Omega_i(t^+)$ knowing $\Omega_i(t^-)$ is given by $k(\bx,\bOmega_i(t^-),\bOmega)d\bOmega$
\end{enumerate}

It can be proved that such a strategy gives a good approximation of the solution $u(t,\bx,\bOmega)$ to
\eqref{E_transport} in the following sense \cite[Theorem 3.2.1]{lapeyre-pardoux-sentis}: assume that 
\begin{equation}\label{eq:mu}
  \eta_N(t,d\bx,d\bOmega) = \alpha \sum_{i=1}^N w_i(t) \delta_{X_i(t),\Omega_i(t)}(d\bx,d\bOmega),
\end{equation}
then this measure converges narrowly to $u(t,\bx,\bOmega)d\bx d\bOmega,$ as $N\to+\infty$. 

\medskip

It remains to include the influence of the source $Q(x)$. For this purpose, we split the time interval into time
steps of equal size $\Delta t$. What follows can easily be generalized to non-constant time steps, but this not our
purpose here. At each time step, we generate $M$ more realizations of another jump Markov process, independently of
the initial ones, as follows:
\begin{enumerate}
\item At time $m\Delta t$, we draw $M$ independent couples $\left(X_j^m, \Omega_j^m\right)$ according to the law
  $\overline Q(\bx)d\bx d\bOmega$.\\
  The weight of each particle is $w_j^m(m\Delta t) = \frac {\Delta t} M$
\item Each of these random variables follow the same evolution as in the preceding case, with positions,
  velocities, weights $X_j^m(t),\Omega_j^m(t),w_j^m(t)$, respectively.
\end{enumerate}

Finally, the measure $\eta_N$ defined by \eqref{eq:mu} is replaced by (here, we assume that $n\Delta t \leq t <
(n+1)\Delta t$)
\begin{equation}
  \label{eq:muQ}
  \eta_{N,\Delta t}(t,d\bx,d\bOmega) = \alpha \sum_{i=1}^N w_i(t) \delta_{X_i(t),\Omega_i(t)}(d\bx,d\bOmega) +
  \beta \sum_{m=0}^n \sum_{j=1}^M w_j^m(t) \delta_{X_j^m(t),\Omega_j^m(t)}(d\bx,d\Omega).
\end{equation}
Here again, this measure narrowly converges to $u(t,\bx,\bOmega)d\bx d\bOmega$, as $n,N,M\to+\infty$, (with
$n\Delta t\to t$) where $u$ solves
\eqref{E_transport}, according to \cite{lapeyre-pardoux-sentis}. 

As we already pointed out, the time step $\Delta t$ may be non-uniform, and, moreover, the number $M$ of particles
generated at each time step may depend on the time step $m$. 

\bigskip

Finally, we point out that we did not take care about boundary conditions. The boundary conditions we aim at
imposing are either free boundary condition, or imposed incoming flux. In the first case, if a process hits the
boundary, it simply vanishes. In the second case, we write the boundary condition as a source consisting of a Dirac mass supported by the boundary,
therefore including it into the source $Q$. 

\subsection{Importance sampling}
\label{sec:importance-sampling}

The method of importance sampling is widely used in many applications of Monte Carlo methods, and in particular in
the case of the simulation of transport equations (see for instance
\cite{alcouffe-85,lux-koblinger,spanier-gelbard}). A nice account of this method in the present context can also be
found in \cite{depinay-these}.

The idea of importance sampling is to introduce an importance function, which we call $I(t,\bx,\bOmega)$, and which
is assumed to be positive. Instead of applying a Monte Carlo method to compute $u$, we are going to apply it to the
function
\begin{displaymath}
  \tilde u(t,\bx,\bOmega) = I(t,\bx,\bOmega) u(t,\bx,\bOmega).
\end{displaymath}
A simple computation gives the equation satisfied by $\tilde u$:
\begin{equation}
  \label{eq:biaisage1}
  \partial_t \tilde u + \bOmega\cdot\nabla \tilde u + \tilde\kappa_t(t,\bx,\bOmega)\tilde u =  \int_{S^2}
  \tilde\kappa_s(t,\bx,\bOmega') \tilde u (t,\bx,\bOmega')\tilde k(t,\bx,\bOmega',\bOmega)d\bOmega' + \tilde
  Q(t,\bx,\bOmega),
\end{equation}
with
\begin{equation}
  \label{eq:kappa_t_modifie}
  \tilde\kappa_t(t,\bx,\bOmega) = \kappa_t - \frac{\partial_t I(t,\bx,\bOmega)}{I(t,\bx,\bOmega)} -
  \frac{\bOmega\cdot\nabla I(t,\bx,\bOmega)}{I(t,\bx,\bOmega)},
\end{equation}
\begin{equation}
  \label{eq:kappa_s_modifie}
  \tilde \kappa_s(t,\bx,\bOmega) = \frac{\kappa_s}{I(t,\bx,\bOmega)} \int_{S^2}k(\bx,\bOmega,\bOmega'')I(t,\bx,\bOmega'')d\bOmega'',
\end{equation}
\begin{equation}
  \label{eq:k_modifie}
  \tilde k(t,\bx,\bOmega',\bOmega) = \left( \int_{S^2} k(\bx,\bOmega',\bOmega'')I(t,\bx,\bOmega'') d\bOmega'' \right)^{-1} k(\bx,\bOmega',\bOmega)I(t,\bx,\bOmega),
\end{equation}
and
\begin{equation}
  \label{eq:source_modifiee}
  \tilde Q(t,\bx,\bOmega) = I(t,\bx,\bOmega) Q(\bx).
\end{equation}

Equation \eqref{eq:biaisage1} has a similar structure as \eqref{E_transport}, and a Monte Carlo method can easily
be designed to compute an approximation of its solution. 
It should be noted that, although the coefficients $\kappa_t$ and $\kappa_s$ where assumed to be constant in
\eqref{E_transport}, the new coefficients defined by \eqref{eq:kappa_t_modifie} and \eqref{eq:kappa_s_modifie} do
depend on $t$, $\bx$ and $\bOmega$ since $I$ does. Similarly, $\tilde k$ depends on $t$ although $k$ does not, and
$\tilde Q$ is now a function of $t,\bx,\bOmega$. This is not a problem for Monte Carlo simulations, the only point
is that equation \eqref{eq:amortissement} should be modified as follows:
\begin{equation}
  \label{eq:amortissement_modifie}
   \tilde w_i(t) = \tilde w_i(t_0) \exp\left( -\int_{t_0}^t \left(\tilde\kappa_t(s,X_i(s),\Omega_i) - \tilde\kappa_s(s,X_i(s),\Omega_i) \right)ds\right). 
\end{equation}
\red{Note that the definition \eqref{eq:kappa_t_modifie} of $\tilde \kappa_t$ does not imply that
  $\tilde \kappa_t\geq 0$. This may be a problem when dealing with Monte Carlo simulations of
  \eqref{eq:biaisage1}. However, if $I$ is solution to the
  adjoint equation (see \eqref{eq:adjointe} below), then $\tilde\kappa_t = \tilde\kappa_s \geq 0$. In the present work, $I$ is not
an exact solution of \eqref{eq:adjointe}, but, as it is pointed out in Remark~\ref{rk:kappa} below, we still have
$\tilde \kappa_t\geq0$ and $\tilde\kappa_s\geq 0$.}
The importance function $I$ should be chosen in such a way that the variance of the computed approximation of
$\tilde u$ has a smaller variance than the approximation of $u$ computed with the method described above.

\subsection{Adjoint equation}
\label{sec:adjoint-equation}

It is known (see \cite{lux-koblinger, spanier-gelbard}) that, in order to have zero variance, the importance function should be
solution to the adjoint equation:
\begin{equation}
  \label{eq:adjointe}
  -\partial_t I -\bOmega \nabla I + \kappa_t I = \kappa_s \int_{S^2} I(t,\bx,\bOmega') k(x,\bOmega,\bOmega')d\bOmega'.
\end{equation}
A rigorous proof of the above fact may be found in \cite{spanier-gelbard}, but let us give a simple argument which
indicates that this is indeed the case. We assume that the spatial domain ${\cal D}$ is a ring between $R_0$ and
$R_1>0$:
\begin{displaymath}
  {\mathcal D} = \left\{ \bx \in \RR^3, \quad R_0< |\bx| < R_1\right\}.
\end{displaymath}
We consider equation \eqref{E_transport} in this domain, with initial condition $g=0$ and boundary conditions 
\begin{displaymath}
  u(t,\bx,\bOmega) =
  \begin{dcases}
    1 &\text{if } |\bx| = R_1, \quad \bOmega\cdot n(\bx) <0, \\
    0 &\text{if } |\bx| = R_0, \quad \bOmega\cdot n(\bx) <0.
  \end{dcases}
\end{displaymath}
Here, $n(\bx)$ is the outer normal unit vector to the domain ${\mathcal D}$ at point $\bx$. Actually, for the
domain we are studying, $n(\bx) = \frac{\bx}{|\bx|}$ if $|\bx|=R_1$ and $n(\bx) = -\frac{\bx}{|\bx|}$ if $|\bx|=R_0$.
We assume that we are interested in computing the flux on the target $|\bx| = R_0$, that is,
\begin{equation}
  \label{eq:flux}
  F(T) = \int_0^T \int_{|\bx|=R_0} \int_{\bOmega\cdot n(\bx)>0} |\bOmega \cdot n(\bx)| u(t,\bx,\bOmega)
  d\bOmega d\bx dt.
\end{equation}
This can be done using the above Monte Carlo method. An estimator of the quantity $F(T)$ is then given by the
following:
\begin{equation}
  \label{eq:flux_estime}
  F = \sum_{j, \ |X_j|=R_0, \ \Omega_j\cdot n(X_j)>0} w_j,
\end{equation}
provided that, at each time step, particles are created with $|X_j|= R_1$, and $\Omega_j$ drawn according to
Lambert cosine law (see \cite{lapeyre-pardoux-sentis}), with initial weights equal to $\frac{\Delta t}M$, where $M$
is the number of particles created at each time step.

Let us now make precise the equation which the importance function $I$ solves: we assume that \eqref{eq:adjointe} is satisfied, and that
the following boundary conditions are imposed
\begin{equation}\label{eq:bc_importance}
  I(t,\bx,\bOmega) =
  \begin{dcases}
    0 &\text{if } |\bx| = R_1, \quad \bOmega\cdot n(\bx) >0, \\
    1 &\text{if } |\bx| = R_0, \quad \bOmega\cdot n(\bx) >0.
  \end{dcases}
\end{equation}

Now, consider a Monte Carlo method applied to $\tilde u$, that is, to equation \eqref{eq:biaisage1}. Here, $\tilde
Q = 0$, but the boundary condition on $\tilde u$ is different:
\begin{displaymath}
  \tilde u(t,\bx,\bOmega) =
  \begin{dcases}
    I(t,\bx,\bOmega) &\text{if } |\bx| = R_1, \quad \bOmega\cdot n(\bx) <0, \\
    0 &\text{if } |\bx| = R_0, \quad \bOmega\cdot n(\bx) <0.
  \end{dcases}
\end{displaymath}
Hence, the boundary data for $\tilde u$ amounts to sampling the distribution $I(t,\bx,\bOmega)$. In particular, a
good choice for the initial weights in such a case is the following:
\begin{displaymath}
  \tilde w_j(t) = \frac{\Delta t}{M}\int_{|\bx|=R_1}\int_{\bOmega\cdot n(\bx)<0} |\Omega\cdot n(\bx)| I(t,\bx,\bOmega)d\bOmega d\bx,
\end{displaymath}
with velocities drawn according to the law ${\displaystyle \frac{|\Omega\cdot n(\bx)|I(t,\bx,\bOmega)\mathds{1}_{\left\{ \bOmega\cdot n(\bx) <0\right\}}}{\int_{\bOmega\cdot n(\bx) <0} |\Omega\cdot n(\bx)|I(t,\bx,\bOmega)d\bOmega}}$. 

First, we note that, since $I$ solves \eqref{eq:adjointe}, we have $\tilde \kappa_s = \tilde \kappa_t$. Hence, the
weight of a particle does not change between shocks. Second,
we point out that all particles go to the target. In order to see this, we assume that it is not the
case. Then, there exists a particle which exits the computation domain through the outer boundary
$|\bx|=R_1$. Denote by $\bx$ its position and $\bOmega$ its direction when it exits. Then, $\bOmega\cdot n(\bx)>0$, $|\bx| = R_1$, and $u(t,\bx,\bOmega)
I(t,x,\bOmega) \neq 0$. But this is impossible since the boundary condition satisfied by $I$ is $I(t,\bx,\Omega) =
0$ for such values or $\bx$ and $\bOmega$.

Hence, the estimated flux $\tilde F$ is now deterministic:
\begin{equation}
  \label{eq:flux_estime_biaise}
  \tilde F = \sum_{j, \ |X_j|=R_0, \ \Omega_j\cdot n(X_j)<0} \tilde w_j = \int_0^T
  \int_{|\bx|=R_1}\int_{\bOmega\cdot n(\bx)<0} \left|\bOmega\cdot n(\bx) \right|
  I(t,\bx,\bOmega)d\bOmega  d\bx dt.
\end{equation}
Finally, multiplying \eqref{E_transport} by $I$, integrating,
and using \eqref{eq:adjointe}, a simple integration by parts proves that
\begin{displaymath}
  \int_0^T \int_{|\bx|=R_1} \int_{\bOmega\cdot n(\bx)<0} |\Omega\cdot n(\bx)| I(t,\bx,\bOmega)d\bOmega d\bx =
  \int_0^T \int_{|\bx|=R_0} \int_{\bOmega\cdot n(\bx)>0} |\Omega\cdot n(\bx)| u(t,\bx,\bOmega)d\bOmega d\bx = F(T).
\end{displaymath}
Hence, if one is able to compute $I$ solution to \eqref{eq:adjointe} with boundary conditions
\eqref{eq:bc_importance}, and to compute the integral on the right-hand side of \eqref{eq:flux_estime_biaise}, then
we have an exact evaluation of the quantity $F(T)$.

\section{Computation of the importance function}
\label{sec:comp-import-funct}

As it was stated in the previous section, computing the solution to the adjoint equation allows to have an
importance function such that the result of importance sampling computations has zero variance. However, solving this equation is at least as difficult as solving \eqref{E_transport}. Therefore, a tractable
approximation of the solution to this equation should be sought in order to be used as an importance function.
It is in general
not possible to compute it exactly, but we will see that in particular situations, simplified expressions may
\red{be} 
derived which give good approximation of the solution to the adjoint problem. First, in
Section~\ref{sec:analyt-expr}, we review the work of \cite{depinay-these}, in which an analytic expression was
derived for the importance function in dimension one. Then, in Section~\ref{sec:spherical-case}, we extend this
analysis to the spherically symmetric case. In such a case, an analytic expression is no longer valid, but a
numerical solution may be computed with the characteristics method, if $\kappa_t, \kappa_s$ and $k$ do not vary
in space. \red{If they do vary in space, we do not know for now how to generalize the calculations of Section~\ref{sec:spherical-case}. As mentionned in
Section~\ref{sec:conclusion} below, a possible way to address this question is to compute numerically the
solution. Doing so, one should be careful to have a good balance between precision (which allows for a significant
improvement of the variance) and computational cost.}

\subsection{Analytic expression : one-dimensional planar case}
\label{sec:analyt-expr}

In his Phd thesis \cite{depinay-these}, J.-M Depinay developed a method to compute an importance function in case
of stationary transport equation in slab geometry. \red{This Section is not directly related to what we
  do in the spherically symmetric case. It is only a simple example in which an explicit computation of the
  importance function $I$ is available. In Section~\ref{sec:spherical-case}, we generalize (to some extent) this
  approach.} In the planar case, the unknown $u$ of \eqref{E_transport} is
assumed to depend only on one space variable $x$, where $\bx = (x,y,z)$. Moreover, the structure of the equation
implies (see  \cite{chandrasekhar}) that $u$ depends
on $\bOmega$ only through the scalar product $\mu = \bOmega \cdot (1, 0, 0)$. Therefore, using these
notations, Equation \eqref{E_transport} reads
\begin{displaymath}
\partial_t u + \mu \partial_x u + \kappa_t u  =\kappa_s \int_{-1}^1  k(x,\mu',\mu)u(x,\mu') d\mu' + Q(x).
\end{displaymath}
Looking for stationary solutions, and assuming that the source is zero ($Q = 0$), this reduces to
\begin{equation}
\mu \partial_x u + \kappa_t u  =\kappa_s \int_{-1}^1  k(x,\mu',\mu)u(x,\mu') d\mu'.
\label{EqStatPlan}
\end{equation}
Hence, the corresponding adjoint problem reads
\begin{equation}
-\mu\partial_x I + \kappa_t I  =\kappa_s \int_{-1}^1  k(x,\mu,\mu')I(x,\mu') d\mu'.
\label{Eq_adjoint_I}
\end{equation}
If one assumes that $I$ is of the form
\begin{equation}
I(x,\mu) = \exp(Kx)\Phi_K(\mu),
\end{equation}
then \eqref{Eq_adjoint_I} implies that $\Phi_K(\mu)$ satisfies the equation 
\begin{equation}\label{eq:2}
\Phi_K(\mu) = \frac{\kappa_s}{\kappa_t-K\mu} \int_{-1}^1 k(x,\mu,\mu')\phi_K(\mu')d\mu', 
\end{equation}
where the parameter $K$ is chosen such that 
\begin{equation}\label{eq:1}
\int_{-1}^1 \Phi_K(\mu)d\mu = 1. 
\end{equation}
It can be proved that such \eqref{eq:1} always has a unique solution $K$ (see \cite[Proposition 7]{depinay-these}),
if $k\in L^\infty$. Moreover, assuming that $k$ is constant (that is, $k=\frac12$), then \eqref{eq:2} reduces to
\begin{displaymath}
  \Phi_K(\mu) = \frac{1}{2}\frac{\kappa_s}{\kappa_t-K\mu}.
\end{displaymath}
Hence, $K$ is the unique solution of equation \eqref{eq:1}, which reads
\begin{displaymath}
  \frac12 \int_{-1}^1 \frac{\kappa_s}{\kappa_t - K\mu}d\mu = 1.
\end{displaymath}
Let us insert $\tilde{u}(x,v) = u(x,v) I(x,v) $ into Equation~\eqref{EqStatPlan}. We have the following modified equation for $\tilde{u}$
\begin{equation}
 \mu\partial_x \tilde{u} + \tilde{\kappa}_t \tilde{u} = \int_{-1}^1 \tilde{\kappa}_s(x,\mu')\tilde{k}(x,\mu',\mu)\tilde{u}(x,\mu') d\mu'.
\label{Eq_biaisee_u}
\end{equation}
with the modified parameters
\[
\begin{aligned}
\tilde{\kappa}_s &= \kappa_t - K\mu, \\
\tilde{\kappa}_t &= \kappa_t - K\mu, \\
\tilde{k}(\mu',\mu) &= \frac{1}{2}\frac{\Phi_K(\mu)}{\Phi_K(\mu')} \frac{\kappa_s}{\kappa_t-K\mu'}.
\end{aligned}
\]
\red{These expressions are \eqref{eq:kappa_t_modifie}, \eqref{eq:kappa_s_modifie}, \eqref{eq:k_modifie}, adapted to the
particular case of 1D planar (slab) geometry.}

Let us point out two important things here: first, the above importance function $I(x,\mu)$ is a solution to the
adjoint equation, but does not in general satisfy appropriate boundary conditions. Therefore, it might result in
poor variance reduction in some situation. However, the tests in \cite{depinay-these} indicate very good efficiency
of the method for a case in which a detector is placed far away from an emitting source. Second, this kind of
solution is related to those exhibited in \cite{case-zweifel} (see also \cite{zweifel-2012}). Such solutions are
eigenvectors of the (adjoint) transport operator, corresponding to the largest possible eigenvalues. This is why
they play an important role here. 

Finally, although the above derivation is done with a stationary transport equation, an implicit time scheme leads,
at each time step, to solving a stationary transport equation. Therefore, \red{although we do not use an
implicit time scheme to solve our transport equation (this would imply additional difficulties that go beyond the
scope of the present work, see \cite{fleck-cummings}, the review paper \cite{wollaber} and the references therein)} the use of the corresponding adjoint
solution in the time-dependent case may prove efficient. This is the strategy we are going to apply in the
spherically symmetric case. 
\subsection{The spherical case}
\label{sec:spherical-case}
Now, we want to compute an importance function in the case of spherical geometry. In such a case, the unknown $u$
is assumed to depend on $\bx$ only through $r = |\bx|$, which in turn implies that it depends on $\bOmega$ only through
$\mu = \bOmega \cdot \frac{\bx}{|\bx|}$. As a consequence, Equation \eqref{E_transport} becomes \cite{castor,chandrasekhar,mihalas,pomraning}:
\begin{equation}
\partial_t u + \mu\partial_r u + \frac{1-\mu^2}{r}\partial_\mu u + \kappa_t u = \kappa_s\int_{-1}^1  k(r,\mu',\mu)u(t,r,\mu') d\mu' + Q(r).\label{EqSpherique}
\end{equation}
Here, $\mu = \cos \theta$, where $\theta$ is the angle formed by the radial direction $\bx$ and the direction $\bOmega$.
\begin{center}
\begin{tikzpicture}[scale=1]
\draw (2,0) arc (0:90:2) ;
\draw (0,0) -- (2,0) ;
\draw (0,0) -- (0,2) ;
\draw (0,0) node[below]{$0$} ;
\draw [dashed] (0,0) -- (1.414213562,1.41421362) ;
\draw (1.414213562,1.41421362) node[right] {$\bx$};
\draw (1.414213562,1.41421362) -- (2.5,2.5) ;
\draw (1.414213562,1.41421362) node {$+$};
\draw (1.6,1.6) arc (0:90:0.2) node[above right] {$ \theta$};
\draw[->] (1.414213562,1.41421362) -- (1.414213562,2.5) ;
\draw (1.414213562,2.5) node[left] {$\bOmega$};
\end{tikzpicture}
\end{center}
As it has been done in the general case, we consider an importance sampling function $I(t,r,\mu)$ and define
$\tilde{u}(t,r,\mu)=u(t,r,\mu)I(t,r,\mu)$ in Equation (\ref{EqSpherique}). Thus, we have the following equation for
$\tilde u$: 
\begin{equation}
\partial_t \tilde{u}+ \mu\partial_r \tilde{u} + \frac{1-{\mu}^2}{r} \partial_\mu \tilde{u} + \tilde{\kappa}_t\tilde{u} = \int_{-1}^{1} \tilde{\kappa}_s(x,\mu') \tilde{k}(x,\mu',\mu)\tilde{u}(x,\mu') d\mu' + \tilde{Q},
\label{Eq_biaisee}
\end{equation}
with 
\begin{displaymath}
\tilde{Q}(t,r,\mu) = Q(r)I(t,r,\mu),   
\end{displaymath}
\begin{equation}\label{eq:kappa_s_modifie_spherique}
\tilde{\kappa}_s(\mu) = \kappa_s\frac{1}{I(t,r,\mu)} \int_{-1}^{1} I(t,r,\mu'') k(r,\mu,\mu'') d\mu'',  
\end{equation}
\begin{equation}\label{eq:kappa_t_modifie_spherique}
\tilde{\kappa}_t(\mu) = \kappa_t - \left( \partial_t I + \mu\partial_r I + \frac{1-{\mu}^2}{r} \partial_\mu I \right)\frac{1}{I(t,r,\mu)} ,  
\end{equation}
\begin{equation}\label{eq:k_modifie_spherique}
\tilde{k}(\mu',\mu) = \frac{I(t,r,\mu)k(\mu',\mu)}{\displaystyle \int_{-1}^{1} I(t,r,\mu'') k(r,\mu',\mu'') d\mu''}  
\end{equation}

The setting we are going to use is the following: we want to reproduce the geometry of an ICF experiment, with a
good statistical convergence on the boundary of the target (or at the ablation front, which is even better). In
order to do so, we assume that the computation domain is
\begin{equation}
  \label{eq:def_domaine}
  {\cal D} = \left\{ r, \quad R_0\leq r \leq R_1\right\},
\end{equation}
where $R_0$ is the radius of the target (or of the ablation front), and $R_1>R_0$ is the outer boundary of the
domain. In all the following, $R_1$ is assumed to be fixed, whereas $R_0 = R_0(t)$ may be a function of time,
reflecting the dynamics of the implosion. An incoming flux is imposed on the outer boundary, while the quantity we
want to compute is the outgoing flux at $r=R_0$. Therefore, Equation \eqref{EqSpherique} is set with the boundary
conditions:
\begin{equation}
  \label{eq:conditions_bord}
  u(r,\mu) =
  \begin{dcases}
    1 &\text{ if } \quad r=R_1, \ \mu <0, \\
    0 &\text{ if } \quad r=R_0, \ \mu>0.
  \end{dcases}
\end{equation}

Equivalently, equation \eqref{Eq_biaisee} is set with the boundary conditions:
\begin{equation}
  \label{eq:conditions_bord_biaisee}
  \tilde u(r,\mu) =
  \begin{dcases}
    I(r,\mu) & \text{ if }\quad r=R_1, \ \mu <0, \\
    0 & \text{ if }\quad r=R_0, \ \mu>0.
  \end{dcases}
\end{equation}

Finally, we assume that $k$ is constant, although this is not a limitation in our strategy:
\begin{displaymath}
  k(r,\mu',\mu) = \frac12.
\end{displaymath}

\subsubsection{Solution of the adjoint equation}
To find $I$, we are going to solve the following adjoint equation 
\begin{equation}
\displaystyle
-\mu \partial_r I - \frac{1-\mu^2}{r}\partial_\mu I + \kappa_t I = \kappa_s \left\langle I \right\rangle ,
\label{EquationAdjointe}
\end{equation}
where
\begin{equation}
\displaystyle
\left\langle I \right\rangle = \frac{1}{2} \int_{-1}^{1} I(r,\mu)d\mu.
\end{equation}
We use a stationary approximation for $I$, although equation \eqref{EqSpherique} is not stationary.
\red{This proves sufficient in the tests below, but this approximation will need to be assessed in
  the presence of material motion, as we point out in the conclusion below.} The reason for
this is twofold: first, some of the numerical tests we are going to use are in fact stationary, and second, even in
the case of a non-stationary situation, the only dependence on time in the model is that of $R_0$, which is assumed
in fact to be constant in each time step of the simulation. Therefore, at each time step, the stationary importance
function should give a good variance reduction.

In order to have a zero variance on the inner boundary $R_0$, the boundary conditions for $I$ should be the
following:
\begin{equation}\label{eq:4}
\displaystyle
\left\{
\begin{aligned}
I(R_0,\mu) & = 1 \text{\ \ if } \mu <0 \\
I(R_1,\mu) &= 0 \text{\ \ if } \mu > 0,
\end{aligned}
\right.
\end{equation}
Considering $S(r) = \kappa_s\left\langle I \right\rangle$ as a source, Equation \eqref{EquationAdjointe} reads :
\begin{equation}
\label{eq:3}
-\mu \partial_r I - \frac{1-\mu^2}{r}\partial_\mu I + \kappa_t I = S(r).
\end{equation}
with the same boundary conditions. We are going to use the method of characteristics to solve \eqref{eq:4}-\eqref{eq:3}. In order to do so,
we change variables, setting $x = r\mu$ and $y = r\sqrt{1-\mu^2}$. The domain is thus (see Figure~\ref{domain}):
\begin{equation}
\displaystyle
\left\{ (x,y)\in\left[-R_1\,;\,R_1\right]\times\left[0\,;\,R_1 \right] \ |\ {R_0}^2 \leq x^2+y^2 \leq {R_1}^2
\right\} = \left\{ (x,y)\in\mathbb{R}\times\mathbb{R}^+ \ |\  {R_0}^2 \leq x^2+y^2 \leq {R_1}^2\right\}.
\end{equation}
We denote by $J$ the new unknown, that is, $I$ as a function of the new variables:
\begin{displaymath}
  I(r,\mu) = J\left(r\mu,r\sqrt{1-\mu^2}\right), \quad \text{and}\quad J(x,y) = I\left(\sqrt{x^2+y^2}, \frac x
  {\sqrt{x^2+y^2}}\right).
\end{displaymath}
Equation (\ref{eq:3}) becomes 
\begin{equation}
-\partial_x J + \kappa_t J = S\left(\sqrt{x^2+y^2} \right).
\end{equation}
This is equivalent to 
\begin{equation}
-\partial_x\left(J(x,y)e^{-\kappa_t x}\right)  = S\left(\sqrt{x^2+y^2} \right)e^{-\kappa_t x}
\label{EquationAdjointeJ}
\end{equation}
with the boundary conditions
\begin{equation}\label{eq:conditions_bord_J}
\displaystyle
\left\{
\begin{aligned}
& J(x,y) = 1 \text{\ \ if } x < 0, \text{\ \ } x^2+y^2={R_0}^2, \\
& J(x,y) = 0 \text{\ \ if } x > 0, \text{\ \ } x^2+y^2={R_1}^2.
\end{aligned}
\right. 
\end{equation}
To solve (\ref{EquationAdjointeJ})-\eqref{eq:conditions_bord_J}, we split the domain into three different parts
(see Figure \ref{domain}): 
\begin{equation}
\displaystyle
\begin{aligned}
{\mathcal D}_1 & = \left\{ (x,y)\quad |\quad x<0,\quad y<R_0, \quad R_0< \sqrt{x^2+y^2}<R_1 \right\}, \\
{\mathcal D}_2 & = \left\{ (x,y)\quad|\quad x>0,\quad y<R_0, \quad R_0< \sqrt{x^2+y^2}<R_1\right\}, \\
{\mathcal D}_3 & = \left\{ (x,y)\quad|\quad  y>R_0, \quad R_0< \sqrt{x^2+y^2}<R_1\right\}. 
\end{aligned}
\end{equation}
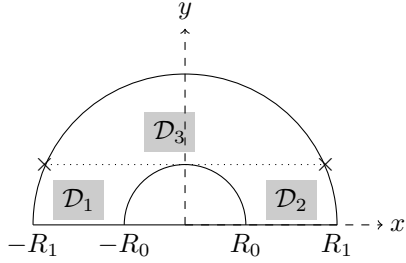
\begin{figure}[!h]
\begin{center}
\begin{tikzpicture}[scale=1]
\draw (0.8,0) arc (0:180:0.8) ;
\draw (2,0) arc (0:180:2) ;
\draw (-2,0) -- (2,0) ;
\draw (-1.85,0.8) node{$\times$} ;
\draw (1.85,0.8) node{$\times$} ;
\draw[dotted] (-1.85,0.8) -- (1.85,0.8);
\draw (0.8,0) node[below]{$R_0$} ;
\draw (-0.8,0) node[below]{$-R_0$};
\draw (2,0) node[below]{$R_1$};
\draw (-2,0) node[below]{$-R_1$};

\draw(-1.4,0.6) node[below, fill=black!20]{${\mathcal D}_1$};
\draw (1.4,0.6) node[below, fill=black!20]{${\mathcal D}_2$};
\draw (-0.2,1.5) node[below,fill=black!20]{${\mathcal D}_3$};

\draw[dashed, ->] (0,0) -- (0,2.6) node[above]{$y$} ; 
\draw[dashed, ->] (0,0) -- (2.6,0) ; 
\draw (2.8,0) node{$x$};
\end{tikzpicture}
\end{center}
\caption{Splitting of the domain in three different zones when applying the method of characteristics.}
\label{domain}
\end{figure} 
In each domain, we integrate the equation along the characteristics and apply the boundary conditions. Thus, we
have, in
domain ${\mathcal D}_1$, 
\begin{displaymath}
J(x,y)= \exp\left( \kappa_t \left( x+\sqrt{{R_0}^2-y^2}\right) \right) + \int_{x}^{-\sqrt{{R_0}^2-y^2}} S\left(\sqrt{s^2+y^2} \right) \exp(\kappa_t (x-s))ds
\end{displaymath}
In domain ${\mathcal D}_2$,
\begin{displaymath}
J(x,y)= \int_{x}^{\sqrt{{R_1}^2-y^2}} S\left(\sqrt{s^2+y^2} \right) \exp(\kappa_t (x-s))ds,
\end{displaymath}
and in domain ${\mathcal D}_3$,
\begin{displaymath}
J(x,y)= \int_{x}^{\sqrt{{R_1}^2-y^2}} S\left(\sqrt{s^2+y^2} \right) \exp(\kappa_t (x-s))ds.
\end{displaymath}
Collecting all these results, the importance function $I(r,\mu)$ reads
\begin{multline}\label{eq:fonction_importance_caracteristiques}
I(r,\mu) =  \mathds{1}_{\left\{ \mu <
    -\sqrt{1-\frac{{R_0}^2}{r^2}} \right\}}\exp\left(\kappa_t\left(r\mu + \sqrt{{R_0}^2-r^2+r^2\mu^2} \right) \right) \\ + \int_{r\mu}^{R(r,\mu)} S\left(\sqrt{s^2+r^2-r^2\mu^2} \right) \exp(\kappa_t (r\mu-s))ds,
\end{multline}
where $R(r,\mu) = - \sqrt{{R_0}^2-r^2+r^2\mu^2}\ \mathds{1}_{\left\{ \mu < -\sqrt{1-\frac{{R_0}^2}{r^2}} \right\}}
+ \sqrt{{R_1}^2-r^2+r^2\mu^2}\ \mathds{1}_{\left\{ \mu > -\sqrt{1-\frac{{R_0}^2}{r^2}} \right\}}$.
\red{In the above formulae (and in the sequel), we use the notation $\mathds{1}$ to indicate the step
  function: for any $m\in\RR$,
  \begin{displaymath}
    \mathds{1}_{\left\{\mu < m\right\}} =
    \begin{dcases}
      1 &\text{ if } \mu < m, \\
      0 & \text{ if } \mu \geq m.
    \end{dcases}
  \end{displaymath}}

\begin{figure}[h]
  \centering
\begin{tikzpicture}[>=triangle 90,scale=0.9]
\draw (1,0) arc (0:180:1) ;
\draw (5,0) arc (0:180:5) ;
\draw (4,0)--(1,0) ;
\draw (5,0)--(4,0)--(-4.5,2.2) arc (154:0:5) -- cycle ;
\draw  [dashed](4,0) -- (2.5,4.33) node[above right] {\scriptsize $R(r,\mu), \ \mu >\mu_d(r)$} ;
\draw (2.5,4.33) node {\scriptsize $\times$}; 
\draw [->,OliveGreen] (4,0)--(2,0.5) node[above] {\scriptsize $\boldsymbol{\mu_d(r)}$};
\draw  [dashed] (4,0) -- (0.866,0.5) node[left,align=center] {\scriptsize $R(r,\mu),$\\ \tiny $\mu<\mu_d(r)$} ;
\draw (0.866,0.5) node {\scriptsize $\times$}; 
\draw (1,0) node {\scriptsize $\times$}; 
\draw (1,0) node[below]{\scriptsize $R_0$} ;

\draw (4,0) node {\scriptsize $\times$}; 
\draw (4,0) node[below]{\scriptsize $r$} ;

\draw (5,0) node {\scriptsize $\times$}; ;
\draw (5,0) node[below]{\scriptsize $R_1$} ;
\end{tikzpicture}  
  \caption{Representation of $R(r,\mu)$ and $\mu_d(r) = -\sqrt{1-\frac{{R_0}^2}{r^2}}$, used in the explicit formula for $I(r,\mu)$~\eqref{eq:fonction_importance_caracteristiques}.}
\end{figure}
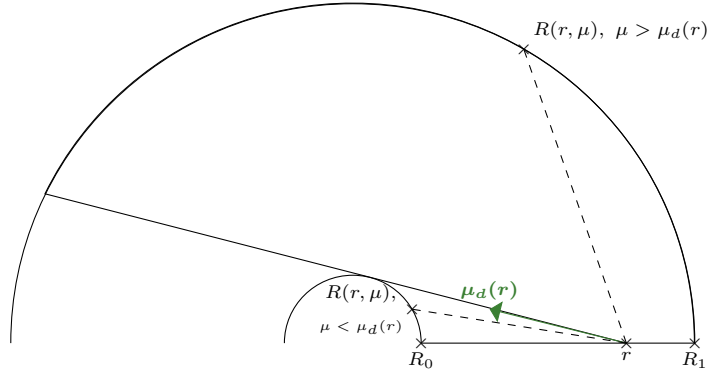

\subsubsection{Integral equation on $S$}
Recalling that ${\displaystyle S(r) = \kappa_s\left\langle I \right\rangle }$, and integrating
\eqref{eq:fonction_importance_caracteristiques} with respect to $\mu$, one obtains an integral equation on $\phi$
defined by
\begin{displaymath}
  \phi(r) = r S(r) = \kappa_s r \langle I\rangle.
\end{displaymath}
This equation reads
\begin{multline}
  \label{eq:equation_integrale}
\phi(r) = \frac{\kappa_s}2 \int_{-1}^{-\sqrt{1-\frac{{R_0}^2}{r^2}}} r\exp\left(\kappa_t\left(r\mu + \sqrt{{R_0}^2-r^2+r^2\mu^2} \right)
\right)d\mu \\+ \frac{\kappa_s}2 \int_{-1}^1 \int_{r\mu}^{R(r,\mu)}\ r\frac{\phi\left(\sqrt{s^2+r^2-r^2\mu^2}
  \right)}{\sqrt{s^2+r^2-r^2\mu^2}} \exp(\kappa_t (r\mu-s))dsd\mu.
\end{multline}
It is possible to compute exactly the integrals with respect to $\mu$, using the exponential integral function (see
\cite{siewert-thomas,erdmann-siewert}, and Appendix~\ref{IntegralEquation} below). This gives
\begin{equation}\label{eq:equation_integrale_2}
\displaystyle
\begin{aligned}
\phi\left( r \right) = \frac{\kappa_s}{4} & \left( \frac{1}{\kappa_t} \bigg[ \exp\left(\kappa_t \theta \right) \bigg]_{-r+R_0}^{-\sqrt{r^2-{R_0}^2}} + \left({R_0}^2 - r^2 \right)\left[ \kappa_t \Ei \left( \kappa_t \theta \right) - \frac{\exp\left( \kappa_t \theta \right)}{\theta} \right]_{-r+R_0}^{-\sqrt{r^2-{R_0}^2}}\right) \\
& + \frac{\kappa_s}{2}\int_{R_0}^{R_1} \phi(r') \bigg[ \Ei \left( \kappa_t \left(-\sqrt{r^2-{R_0}^2} - \sqrt{{r'}^2-{R_0}^2} \right) \right) - \Ei \big( \kappa_t\left(-| r -r' | \right)\big)   \bigg] dr'.
\end{aligned}
\end{equation}
Here, $\Ei(x)$ is the exponential integral defined by
\begin{equation}
\displaystyle
\Ei \left(x\right) = \int_{-x}^{\infty} \frac{\exp\left( -t\right)}{t}dt.
\end{equation}

\subsubsection{Numerical computation of $I$}

\emph{A priori,} it is not possible to solve this integral equation exactly. However, it is possible to solve it numerically. For this purpose, we introduce a mesh to discretize the space variable $r$. In all the numerical
examples we are going to give, we use a uniform mesh, although this is not essential. Thus, we use the following
notations: let $N_r$ be a positive integer, and define $\Delta r = \frac{R_1-R_0}{N_r}$. Note however that in some
of the cases treated below, $R_0$ depends on $t$. In such cases, we use a mesh independent of $t$, with $\Delta r = \frac{R_1}{N_r}$.

\begin{multline}
  \label{eq:notation_mesh}
  \forall\ 0\leq j \leq N_r-1, \quad M_j = \left[r_{j-1/2},r_{j+1/2}\right], \quad r_{j-1/2} = R_0 + j\Delta r, \quad
  r_{j+1/2} = R_0 +  (j+1)\Delta r, \\ r_j = R_0+\left(j+\frac12\right)\Delta r.
\end{multline}
We then use a piecewise constant approximation of $\phi$, defining
\begin{displaymath}
  \phi(r) = \sum_{j=0}^{N_r-1} \phi_j \mathds{1}_{M_j}(r), \quad \text{hence} \quad \int_{M_j} \phi = |M_j|\phi_j =
  \Delta r \phi_j.
\end{displaymath}
Inserting this into \eqref{eq:equation_integrale_2}, and using a piecewise constant approximation for all the
functions appearing in the integrals, we infer
\begin{displaymath}
  \phi_j = b_j + \frac{\kappa_s}2 \sum_{i=0}^{N_r-1} \phi_i \Delta r\left[ \Ei\left(-\kappa_t \sqrt{r_j^2 - R_0^2} -
      \kappa_t\sqrt{r_i^2 - R_0^2}\right) - \Ei\left(-\kappa_t \left|r_j-r_i\right| \right)\right],
\end{displaymath}
where
\begin{equation}\label{eq:second_membre}
  b_j = \frac{\kappa_s}{4} \left( \frac{1}{\kappa_t} \bigg[ \exp\left(\kappa_t \theta \right) \bigg]_{-r_j+R_0}^{-\sqrt{r_j^2-{R_0}^2}} + \left({R_0}^2 - r_j^2 \right)\left[ \kappa_t \Ei \left( \kappa_t \theta \right) - \frac{\exp\left( \kappa_t \theta \right)}{\theta} \right]_{-r_j+R_0}^{-\sqrt{r_j^2-{R_0}^2}}\right)
\end{equation}
Hence, we are lead to the following linear system satisfied by $\phi$:
\begin{equation}\label{eq:systeme_lineaire_phi}
\left(I_d - A \right)\phi = b,
\end{equation} 
where $I_d$ is the identity matrix and $A$ is defined by:
\begin{equation}\label{eq:matrice_A}
A_{ij} = \frac{\kappa_s}2\Delta r\left[
 \Ei \left( \kappa_t \left(-\sqrt{{r_i}^2-{R_0}^2} - \sqrt{{r_j}^2-{R_0}^2} \right) \right) - \Ei \big( \kappa_t\left(-| r_i -r_j | \right)\big)\right],
\end{equation} 
and the right-hand side $b$ is defined by \eqref{eq:second_membre}. This formula is valid only in the case $i\neq
j$. If $i=j$, the singularity of $\Ei$ at the origin does not allow for the use of \eqref{eq:matrice_A}. In order
to compute them, we note that $H=1$ is the unique solution of the system
\[
\left\{
\begin{aligned}
& -\mu \partial_r H - \frac{1-\mu^2}{r}\partial_\mu H + \kappa_t H = \kappa_s \left\langle H \right\rangle + \kappa_t-\kappa_s \\
&H\left(R_0,\mu \right)  = 1 \text{ , } \mu<0\\
&H\left(R_1,\mu \right)  = 1 \text{ , } \mu>0\\
\end{aligned}
\right.
\] 
Now, applying the method of characteristics as above to this system, we have an equation for $H$ similar to
\eqref{eq:fonction_importance_caracteristiques}, with an additional term due to the boundary condition at $r=R_0$:
\begin{multline} \displaystyle
H\left(r,\mu\right) = \exp\left( \kappa_t\left( r\mu + \sqrt{ r^2\mu^2 - r^2 + {R_0}^2}
  \right)\right)\mathds{1}_{\left\{ \mu < -\sqrt{1-\frac{{R_0}^2}{r^2}}\right\}} \\
+  \exp\left( \kappa_t\left( r\mu + \sqrt{ r^2\mu^2 - r^2 + {R_1}^2} \right)\right)\mathds{1}_{\left\{ \mu > -\sqrt{1-\frac{{R_0}^2}{r^2}}\right\}} \\ + \int_{r\mu}^{R(r,\mu)} \left(\frac{\Psi\left(\sqrt{s^2+r^2-r^2\mu^2}\right)}{\sqrt{s^2+r^2-r^2\mu^2}}+ \kappa_t-\kappa_s\right)\exp\left(\kappa_t \left(r\mu-s\right) \right)ds.
\end{multline}
Moreover, integrating with respect to $\mu$, we also have a relation similar to \eqref{eq:equation_integrale_2}:
\begin{multline}
\label{eq:equation_integrale_Psi} \displaystyle
\Theta\left( r \right):=\kappa_sr\langle H\rangle = \frac{\kappa_s}{4}  \left( \frac{1}{\kappa_t} \bigg[ \exp\left(\kappa_t \theta \right) \bigg]_{-r+R_0}^{-\sqrt{r^2-{R_0}^2}} + \left({R_0}^2 - r^2 \right)\left[ \kappa_t \Ei \left( \kappa_t \theta \right) - \frac{\exp\left( \kappa_t \theta \right)}{\theta} \right]_{-r+R_0}^{-\sqrt{r^2-{R_0}^2}}\right) \\
+ \frac{\kappa_s}{4}  \left( \frac{1}{\kappa_t} \bigg[ \exp\left(\kappa_t \theta \right) \bigg]^{r+R_1}_{-\sqrt{r^2-{R_0}^2}+\sqrt{{R_1}^2-{R_0}^2}} + \left({R_1}^2 - r^2 \right)\left[ \kappa_t \Ei \left( \kappa_t \theta \right) - \frac{\exp\left( \kappa_t \theta \right)}{\theta} \right]^{r+R_1}_{-\sqrt{r^2-{R_0}^2}+\sqrt{{R_1}^2-{R_0}^2}}\right) \\
+ \frac{\kappa_s}{2}\int_{R_0}^{R_1} \bigg(\Psi(r') + \left(\kappa_t - \kappa_s \right)r'\bigg) \bigg[ \Ei \left( \kappa_t \left(-\sqrt{r^2-{R_0}^2} - \sqrt{{r'}^2-{R_0}^2} \right) \right) - \Ei \big( \kappa_t\left(-| r -r' | \right)\big)   \bigg] dr'
\end{multline}
Using the fact that $H=1$, and assuming a piecewise constant approximation of $\Theta$, we may assume
\begin{equation}\label{eq:4.1}
  \Theta_j = \kappa_s r_j.
\end{equation}
Inserting this into \eqref{eq:equation_integrale_Psi}, we find
\begin{displaymath}
\kappa_s\left( r_j - \sum_{i=0}^{N_r-1}A_{ji}r_i\right) = b_j + c_j + \sum_{i=0}^{N_r-1} A_{ji}d_i\ ,
\end{displaymath}
where the coefficients $A_{ij}$ and $b_j$ are defined as above by \eqref{eq:second_membre} and
\eqref{eq:matrice_A}. Here, the coefficients $c_j$ are given by 
\begin{displaymath}
c_j  = \frac{\kappa_s}{4}  \left( \frac{1}{\kappa_t} \bigg[ \exp\left(\kappa_t \theta \right)
  \bigg]^{r_j+R_1}_{-\sqrt{r_j^2-{R_0}^2}+\sqrt{{R_1}^2-{R_0}^2}} + \left({R_1}^2 - r_j^2 \right)\left[ \kappa_t
    \Ei \left( \kappa_t \theta \right) - \frac{\exp\left( \kappa_t \theta \right)}{\theta}
  \right]^{r_j+R_1}_{-\sqrt{r_j^2-{R_0}^2}+\sqrt{{R_1}^2-{R_0}^2}}\right),
\end{displaymath}
and 
\begin{displaymath}
d_j  = \left(\kappa_t - \kappa_s \right)r_j.
\end{displaymath}
Hence, the vector $\left(\Theta_j\right)$ satisfies the equation
\begin{displaymath}
\left(I_d - A \right)\Theta = b + c + Ad,
\end{displaymath}
Hence, using \eqref{eq:4.1},
\begin{displaymath}
  r_j\kappa_s = b_j + c_j + \kappa_t \sum_{i=0}^{N_r-1} A_{ji}r_i.
\end{displaymath}
Finally,
\begin{equation}\label{eq:11}
A_{jj} = \frac{r_j\kappa_s - b_j - c_j}{\kappa_t r_j} - \sum_{i\neq j} A_{ji}\frac{r_i}{r_j}.
\end{equation}
We have an expression of the diagonal coefficients $A_{jj}$ in terms of the off-diagonal ones.

\medskip

After solving \eqref{eq:systeme_lineaire_phi},
we use $\phi$ to define an approximation of the importance function $I(r,\mu)$:
\begin{multline}\label{eq:7}
I\left(r,\mu\right) = \mathds{1}_{\left\{ \mu < -\sqrt{1-\frac{{R_0}^2}{r^2}}\right\}} \exp\left( \kappa_t\left( r\mu + \sqrt{ r^2\mu^2 - r^2 + {R_0}^2} \right)\right) \\ +  \sum_{i=0}^{N_r-1} \phi_i\int_{r\mu}^{R(r,\mu)} \mathds{1}_{M_i}\left(\sqrt{s^2+r^2-r^2\mu^2} \right) \frac{\exp\left(\kappa_t \left(r\mu-s\right) \right)}{\sqrt{s^2+r^2-r^2\mu^2}}ds.
\end{multline}

\bigskip

Since $\tilde\kappa_s$ depends on $\mu$, we need a mesh in propagation direction $\mu$. We use a piecewise constant approximation
of $I$ on each cell (both in space and direction). This implies discontinuities at the boundaries of the cells, so
we have to change the value of $\tilde u$ when going from a cell to its neighbour. Indeed, between a cell $M_1$ and
a cell $M_2$, the continuity of $u$ implies that 
\begin{displaymath}
\frac{\tilde u_{M_1}(r,\mu)}{I_{M_1}(r,\mu)} = \frac{\tilde u_{M_2}(r,\mu)}{I_{M_2}(r,\mu)}.
\end{displaymath}
In order to take this into account, we multiply the weight of a particle going from $M_1$ to $M_2$ by
$\frac{I_{M_2}(r,\mu)}{I_{M_1}(r,\mu)}$. 
\red{
\begin{remarque}
  \label{rk:kappa}
A simple computation proves that, if $I$ is an exact solution of \eqref{EquationAdjointe}, then $\tilde\kappa_t =
\tilde\kappa_s \geq 0$. Here, we do not have an exact solution, but one can still prove, with the same computation,
that $\tilde\kappa_s \geq 0$ and $\tilde\kappa_t \geq 0$. 
\end{remarque}}
\red{
\noindent\emph{Proof:} If $I$ is solution to \eqref{EquationAdjointe}, then $I\geq 0$. Hence, by definition, we have $\tilde\kappa_s
\geq 0$, according to \eqref{eq:kappa_s_modifie_spherique}. Now, considering $\tilde \kappa_t$, we have, using
\eqref{eq:kappa_t_modifie_spherique} and \eqref{EquationAdjointe},
\begin{displaymath}
  \tilde \kappa_t = \kappa_t -\frac 1 I \left(\mu\partial_r I + \frac {1-\mu^2}r \partial_\mu I\right) = \kappa_t -
  \frac 1 I \left( \kappa_t I - \kappa_s \langle I \rangle\right) = \kappa_s\frac{\langle I \rangle}{I} \geq 0.
\end{displaymath}
Next, we consider the case in which \eqref{EquationAdjointe} is replaced by \eqref{eq:3}, where $S$ is no more equal to
$\kappa_s \langle I \rangle$. However, $S(r)=\phi(r)/r$, where $\phi$ is numerically computed by solving
\eqref{eq:systeme_lineaire_phi}, where $A$ is defined by \eqref{eq:matrice_A} and \eqref{eq:11}, and $b$ by
\eqref{eq:second_membre}. With these definitions and the fact that the function $\Ei$ is negative and decreasing on
$\RR^{-}$, we infer that $I_d -A$ is an M-matrix, and that $b_j\geq 0$, for all $j$. Hence, $S\geq 0$, from which we
deduce again that $I\geq 0$. The proof of $\tilde \kappa_s\geq 0$ is exactly the same as above. Turning to $\tilde\kappa_s
\geq 0$, we point out that $S\geq 0$, hence the same computation as above gives
\begin{displaymath}
  \tilde \kappa_t = \kappa_t -\frac 1 I \left(\mu\partial_r I + \frac {1-\mu^2}r \partial_\mu I\right) = \kappa_t -
  \frac 1 I \left( \kappa_t I - S\right) = \frac{S}{I} \geq 0.
\end{displaymath}
\hfill$\Box$
}
\section{Numerical results}
\label{sec:numerical-results}
The law of large numbers states that for independent and identically distributed random variables, the sample
average converges to the expected value when the number of random variables increases. \red{(In the
  present context, a random variable is synonymous to a Monte Carlo particle.)} In addition, the central
limit theorem implies that the rate of convergence is ${\displaystyle \frac{1}{\sqrt{N}} }$, where $N$ is the
number of random variables. In our case, $N$ is equal to the number of particles used in the simulation. 
For Monte Carlo methods, the computational cost is proportional to the number $N$ of random variables. So, to
compare the performance of different methods, we have to take the calculation time $T$ into account. Knowing that
$T\approx Nt $, where $t$ is the time to generate one realization of a random variable, we define the figure of
merit (F.O.M) of a method by 
\begin{equation}\label{eq:FOM}
\text{F.O.M} = \frac{1}{\sigma^2 t} = \frac{1}{\Sigma_N^2T}, 
\end{equation}
with $\sigma^2$ is the variance of $X_i$, where $\left( X_i\right)_{i \in \mathbb{N}}$ is used in the simulation sequence of independent random variables. An unbiased estimator for $\sigma^2$ is 
\[ \displaystyle
\sigma_N^2 = \frac{1}{N-1}\sum_{i=1}^{N} \left(X_i - \bar{X}_N \right)^2 \text{\; with \; }  \displaystyle \bar{X}_N = \frac{1}{N}\sum_{i=1}^{N}X_i. 
\]
Moreover $\Sigma_N^2 $ is the variance of $\bar{X}_N$. An unbiased estimator for $\Sigma_N^2$ is $\displaystyle
\frac{\sigma_N^2}{N}$. \red{In steady cases below, we apply the above formulas. In unsteady cases, we
  apply them to time-integrated values.}

\bigskip

We are now going to give some numerical results obtained with the method developed so far. First, we provide two
\red{verification} cases, which indicate that our implementation of the Monte Carlo method is correct. In these cases, we
have analytical solutions, allowing to assess the statistical convergence of the method. Second, we provide
variance reduction tests, in which we do not have any analytical solution, so we only study variance reduction when
importance sampling is applied. The first case is stationary, and the second one is unsteady, with data in
agreement with FCI simulations. 

\subsection{A stationary \red{verification} test case}
\label{sec:Cas_Test_Siewert_Thomas}

This test is borrowed from \cite{siewert-thomas} and \cite{erdmann-siewert}, and is used as a \red{verification} procedure for
our Monte Carlo code (without using the importance sampling method). We solve Equation \eqref{EqSpherique} with a
point source located at $r=R_{\rm source}$ \red{(actually, this is a point source only if $R_0=0$, but we
  nevertheless use this denomination even if $R_0>0$.)}:
\[
Q(r) = \frac{1}{4\pi r^2}\delta\left(r-R_{\rm source} \right),
\]
The domain is $(0,1)$, that is, \eqref{eq:def_domaine} with $R_0=0$ and $R_1=1$. The cross sections are such that
$\kappa_t = 1$. We test $\kappa_s = 0.3$ and $\kappa_s = 0.9.$ The boundary conditions correspond to zero incoming
flux:
\begin{displaymath}
  u(R_1,\mu) = 0, \quad \forall \mu<0. 
\end{displaymath}
We use $10^5$ particles, and compute the zero-moment of the intensity $\psi$ as a function of $r$, on a mesh with $1000$ identical
cells, that is, $\Delta r= 10^{-3}$. Recall that
\begin{equation}\label{eq:6}
  \psi(r) = 2\pi \int_{-1}^1 u(r,\mu)d\mu,
\end{equation}
so that an unbiased estimator of the average value of $\psi$ on the cell $M$ is given by 
\begin{displaymath}
  \psi_j = \sum_{X_j\in M} w_j \ \mathop{\longrightarrow}_{N\to+\infty}\  \int_M \psi,
\end{displaymath}
where $j$ is the index of a particle, $X_j$ its position and $w_j$ its weight. 

\begin{figure}[!h]
   \centering
   \subfigure[$R_{\rm source}=0{.}05$]{%
      \includegraphics[scale=0.30]{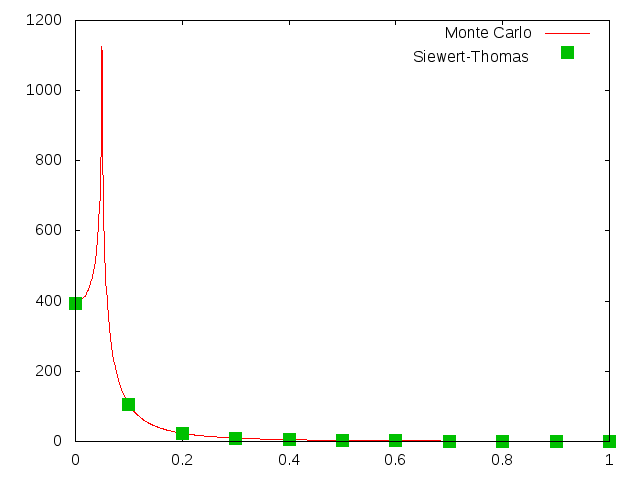}}%
   \hspace{1mm}%
   \subfigure[$R_{\rm source}=0{.}45$]{%
      \includegraphics[scale=0.30]{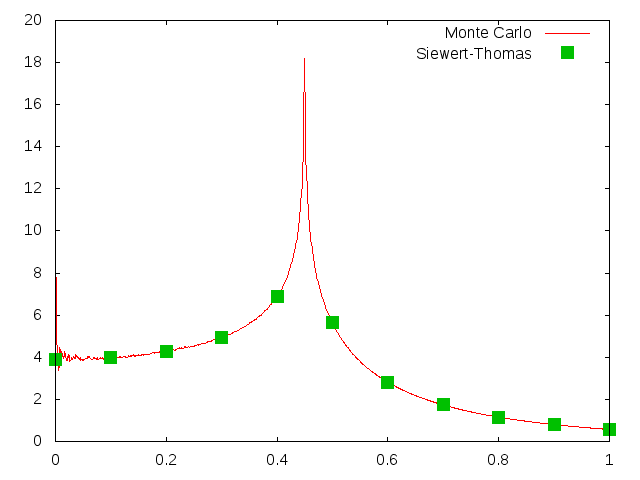}}%
  \subfigure[$R_{\rm source}=0{.}95$]{%
      \includegraphics[scale=0.30]{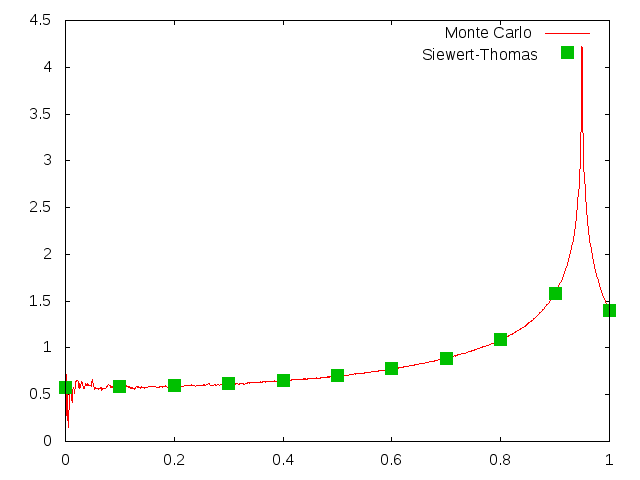}}%
   \hspace{1mm}%
\caption{\label{Siewert-ThomasC1}The flux $\psi$ defined by \eqref{eq:6}: comparison between the analytical
  solution of Siewert and Thomas \cite{siewert-thomas}, and the result of our Monte Carlo code \red{without
    importance sampling}. Here, $\kappa_s=0.3$ and $R_1=1$.}
\end{figure}

\begin{figure}[!h]
   \centering
   \subfigure[$R_{\rm source}=0{.}05$]{%
      \includegraphics[scale=0.30]{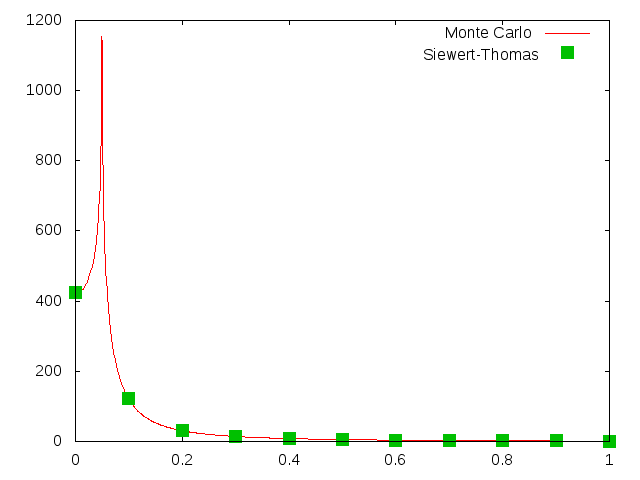}}%
   \hspace{1mm}%
   \subfigure[$R_{\rm source}=0{.}45$]{%
      \includegraphics[scale=0.30]{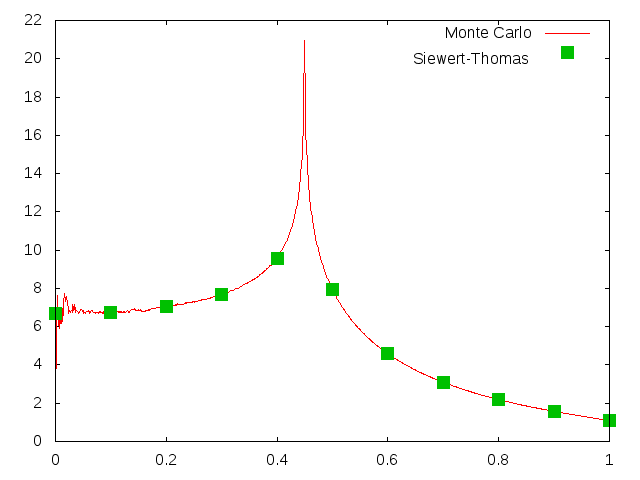}}%
  \subfigure[$R_{\rm source}=0{.}95$]{%
      \includegraphics[scale=0.30]{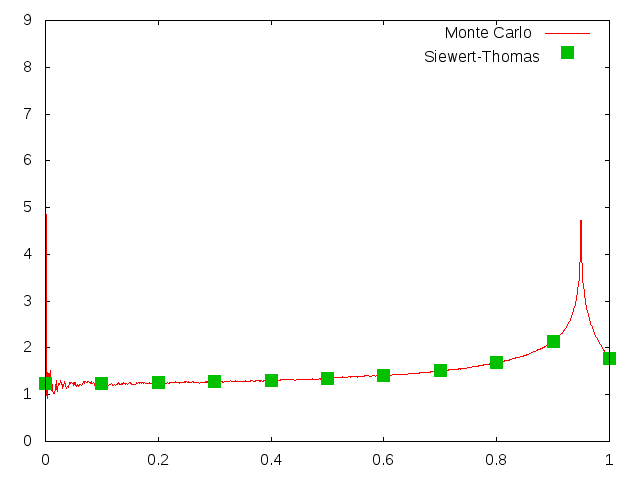}}%
   \hspace{1mm}%
   \caption{\label{Siewert-ThomasC2}The flux $\psi$ defined by \eqref{eq:6}: comparison between the analytical
  solution of Siewert and Thomas \cite{siewert-thomas}, and the result of our Monte Carlo code \red{without
    importance sampling}. Here, $\kappa_s=0.9$ and $R_1=1$.}
\end{figure}

A semi-analytical solution is derived in \cite{erdmann-siewert}. This solution involves an integral which is
computed numerically. In Figure~\ref{Siewert-ThomasC1} and Figure~\ref{Siewert-ThomasC2}, we compare the result of
our Monte-Carlo code \red{(without importance sampling)} with this analytical solution. 

This test shows a good agreement between the result of our code and the analytical solution. The statistical noise
is more important in the cells near the origin. This can be explained by the fact that these cells are small, so
very few particles are present in them. The singularity at $r=R_{\rm source}$ is well reproduced.

\subsection{An unsteady \red{verification} test case}
\label{sec:Cas_Test_Validation_Instationnaire}

In this test, we assume that $R_1 = 1$, and that the target is locate at $r=R_0(t)$, where we have set
\begin{equation}\label{eq:10}
  R_0(t) = \alpha + \beta t, \quad \alpha = 0.37625, \quad \beta = -0.027625
\end{equation}
These data are borrowed from physically relevant cases of inertial confinement fusion (after
adimensionalization). \red{The incoming flux imposed at $r=R_1$ is equal to $1$ (with Lambert cosine law) between
  times $0$ and $T_{\rm max}$, and $0$ between $T_{\rm max}$ and $T=10$. We tune $T_{\rm max}$ so as to have an
  exact value for the flux on the target integrated between $0$ and $T$, as is explained below.}

\begin{figure}[ht]
\centering
	\includegraphics[scale=0.5]{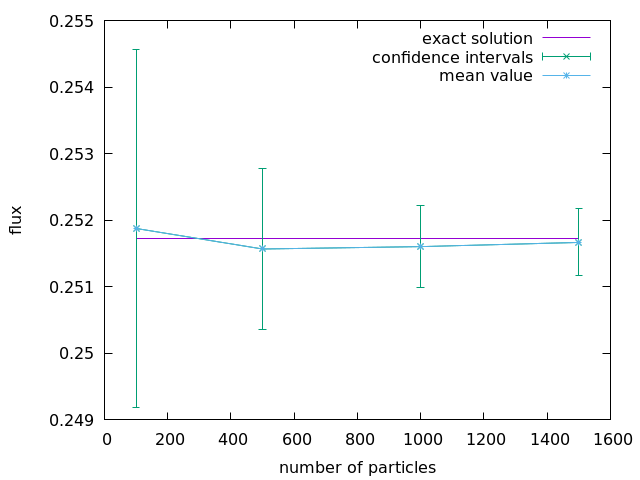}
	\caption{\red{Verification} test case for the Monte Carlo code \red{(without importance sampling)} with $\kappa_s=0$ and $\kappa_t=0$}
	\label{Validation_Monte_Carlo_S0}
\end{figure}

In the case $\kappa_s=0$ and $\kappa_t=0$, we have an analytic expression for the solution, and the flux at the
boundary of the target, integrated in time from $0$ to $T_{\rm max}$, is equal to
\begin{equation}\label{eq:5}
F = \int_{0}^{T_{\rm max}} \int_{-1}^{\mu(t)} |\mu| d\mu dt,
\end{equation}
where the time $T_{\rm max}$ is made precise below. 
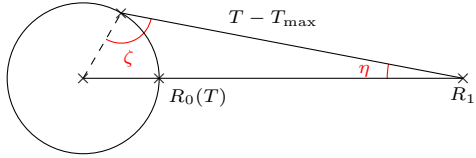
\begin{figure}[h!t]
  \centering
\begin{tikzpicture}[>=triangle 90,scale=1]
\draw [color=red] (4,0) arc (180:169:1cm);
\draw [color=red] (3.7,0.3) node[below]{\scriptsize $\eta$} ;
\draw [color=red] (0.3,0.5196) arc (240:350:0.4cm);
\draw [color=red] (0.6,0.5) node[below]{\scriptsize $\zeta$} ;
\draw (0,0) circle (1) ;
\draw (0,0) node {\scriptsize $\times$};
\draw (5,0) -- (0.5,0.866025) ;
\draw (5,0) -- (0,0) ;
\draw (0.5,0.866) node {\scriptsize $\times$}; 
\draw (2.5,0.8) node {\scriptsize $T-T_{\rm max}$} ;
\draw  [dashed] (0,0) -- (0.5,0.866025) ;
\draw (1,0) node {\scriptsize $\times$}; 
\draw (1,0) node[below right]{\scriptsize $R_0(T)$} ;
\draw (5,0) node {\scriptsize $\times$}; ;
\draw (5,0) node[below]{\scriptsize $R_1$} ;
\end{tikzpicture}  
  \caption{$T_{\rm max}$ is the time such that particles generated at $t<T_{\rm max}$ with initial direction $\mu
    \in [-1,\mu(t)]$ reach the target. \red{Here, $\cos\zeta = \beta$, where $\beta$ is defined by \eqref{eq:10}, and $\cos\eta = -\mu(t)$.}}
\end{figure}
The direction  $\mu(t)$ is, for a particle generated at time $t$ and reaching the target, the largest possible
propagation direction. In order to compute it, we compute the trajectory of the corresponding particle, defined by $\dot r =
\mu$ and $\dot \mu = \frac{1-\mu^2}{r},$ hence
\begin{displaymath}
  r(t+t')^2 = R_1^2 +2t'R_1\mu(t) +t'^2, \quad \mu(t+t') = \frac{R_1\mu(t)+t'}{\sqrt{R_1^2 +2t'R_1\mu(t) +t'^2}}.
\end{displaymath}
This trajectory crosses the inner boundary if and only if the equation $r(t+t') = R_0(t+t')$ has a solution. This
equation reads $R_0(t+t')^2 =  R_1^2 +2t'R_1\mu(t) +t'^2,$ which is a second-degree equation in $t'$. The maximum
direction $\mu(t)$ corresponds to the case when the discriminant is $0$. Computing it, we find
\begin{equation}\label{eq:8}
\mu(t) = \frac{1}{R_1}\left[\beta R_0(t) - \sqrt{\left(1-\beta^2 \right)\left(R_1^2 - R_0(t)^2 \right)} \right].
\end{equation}
Hence, \eqref{eq:5} also reads
\begin{equation}\label{eq:9}
F = \frac{1}{2}\int_{0}^{T_{\rm max}} \left( 1-\mu(t)^2\right) dt = \frac12 \int_0^{T_{\rm max}} \left( 2
   -\beta^2 -\frac{R_0(t)^2}{R_1^2} -
    2\beta\frac{R_0(t)}{R_1}\sqrt{\left(1-\beta^2 \right)\left(1 - \frac{R_0(t)^2}{R_1^2} \right)} \right)dt 
\end{equation}
In \eqref{eq:9}, the time $T_{\rm max}$ is such that a particle generated at $t<T_{\rm max}$ reaches the target
with initial direction $\mu\in [-1,\mu(t)]$, where $\mu(t)$ is given by~\eqref{eq:8}. Indeed, for such value of $T_{\rm max}<T$, formula \eqref{eq:5} is valid, whereas if $T_{\rm max}$ is
chosen to be larger, some particles, generated between $T_{\rm max}$ and $T$, never reach the target. Therefore the
value of $\mu(t)$ is no more given by \eqref{eq:8}, and formula \eqref{eq:5} should be modified. In order to avoid
technical difficulties associated to this new value of $\mu(t)$, we
restrict the time integral to $[0,T_{\rm max}].$ A simple computation shows that
\begin{displaymath}
T_{\rm max} = T - R_0(T)\beta - \sqrt{R_0(T)^2\beta^2 - R_0(T)^2+R_1^2}.
\end{displaymath}
Using the above values of $\alpha$ and $\beta$, we find $T_{\rm max} = 9.00777122797$, for which we find
\begin{displaymath}
  F = 0.25172763696.
\end{displaymath}
This computation is carried out in the following conditions: $\Delta r = 10^{-2}$, $\Delta t = 10^{-3}$. The
results are displayed in Figure \ref{Validation_Monte_Carlo_S0}, showing statistical convergence to the exact value
as the number of particles grows. Note that the abscissa in Figure~\ref{Validation_Monte_Carlo_S0} is the number of
particles generated at each time step. Therefore, the total number of particles in the simulation is equal to
$N\times T/\Delta t = N\times 10^4$. 

\bigskip

\begin{figure}[ht]
\centering
	\includegraphics[scale=0.5]{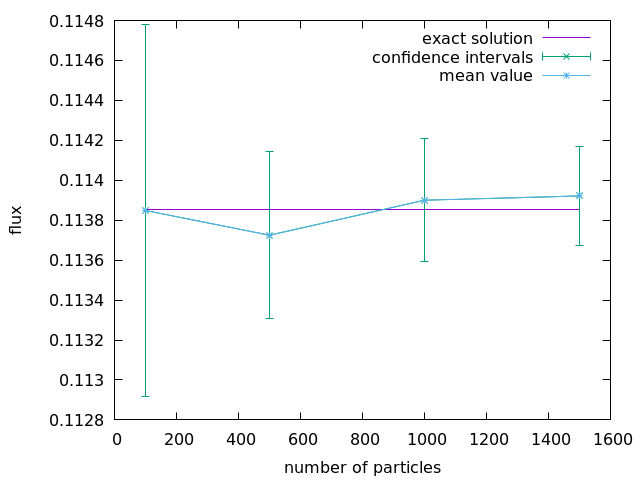}
	\caption{\red{Verification} test for Monte Carlo code \red{(without importance sampling)} with $\kappa_s=0$ and $\kappa_t=1$}
	\label{Validation_Monte_Carlo_S1}
\end{figure}

In the case $\kappa_s=0$ and $\kappa_t=1$, we still have an exact expression of the flux on the inner ball $R_0(t)$:
\begin{displaymath}
F = \int_{0}^{T_{\rm max}} \int_{-1}^{\mu(t)} |\mu|\exp\left[-\left(\kappa_t -\kappa_s\right)\left(\tau(t,\mu) -t \right)\right] d\mu dt
\end{displaymath}
with
\begin{displaymath}
\tau (t,\mu) = \frac{1}{1-\beta^2} \left(t+\alpha\beta - R_1\mu-\sqrt{\left(R_1\mu-t-\alpha\beta \right)^2 - \left( 1-\beta^2 \right)\left(t^2+R_1^2-\alpha^2-2R_1\mu t \right)}\right).
\end{displaymath}
This integral is not explicit, so we applied a numerical integration method to compute it. The result is
\begin{displaymath}
  F =  0.11385526445,
\end{displaymath}
up to an estimated error of $10^{-6}$. The parameters for this test are as follows: $100$ cells in $\Delta r =
10^{-2}$,  $\Delta t= 10^{-3}$. 
Figure~\ref{Validation_Monte_Carlo_S1} gives the results of this test, showing statistical convergence as the
number of particles grows. 

\subsection{Stationary test case: variance reduction}
\label{sec:Cas_Test_Stationnaire}

We now consider a stationary test case again. Here, we solve the transport equation in the domain $\left(R_0,R_1\right)$, with
$R_0 = 0.1$ fixed, and $R_1 = 1$. An incoming flux is imposed on the outer boundary $r=R_1$, of value $1$. The
discretization corresponds to $\Delta r = 10^{-2}$ and $\Delta \mu = 2\times 10^{-3}$. Note that the discretization
in $\mu$ is only used when the importance sampling method is applied.

\subsubsection{Case $\kappa_s=0.9$, $\kappa_t=1$}

\begin{table}[!ht]
  \begin{center}
	\begin{tabular}{|c|c|c|c|c|c|c|}
      \hline
      N & Flux & Variance & Standard deviation & Time & F.O.M & P(N)\\  
      \hline
      100 & 0.038368316 & 1.772142E-03 & 0.04209682 &  0.1731 & 3.26E+04 & 0.87\% \\
      \hline
      500 & 0.038080894 & 2.832987E-04 & 0.01683148 &  0.8360  & 4.22E+04 & 0.72\% \\
      \hline
      1000 & 0.036268938 & 1.554811E-04 & 0.012469206 & 1.3283 & 4.84e+04 &0.95\% \\
      \hline
      1500 & 0.036145858 & 6.882991E-05 &0.00829638 & 1.5196  & 9.56E+04 & 0.82\% \\
      \hline
      10000& 0.040297305 & 1.147480E-05 &0.00338745 & 13.571 & 6.42E+04 & 0.91\% \\
      \hline
      100000&0.03764841 & 1.6171619E-06 &0.00127168 &  90.81 & 6.81E+04 & 0.85\% \\
      \hline
	\end{tabular}
  \end{center}
  \caption{The case $\kappa_s=0.9$, $\kappa_t = 1$ without importance sampling.}
  \label{Tableau_ResultatStat_MC_Ks09_Kt1}
\end{table}  

\begin{table}[!ht]
	\begin{tabular}{|c|c|c|c|c|c|c|c|c|}
  		\hline
  		N & Flux & Variance & Standard deviation & Time1 & Time2 & F.O.M1 & F.O.M2 & P(N)\\  
  		\hline
100 & 0.0386001164 & 1.064560E-05 & 3.2627599E-03 & 4.4519 &0.7602 & 2.11E+05 & 2.07E+06  & 87.7\% \\
		\hline
500 & 0.0376133306 & 8.991524E-07 & 9.4823647E-04 & 7.7332 & 4.0286 & 1.44E+06 & 2.76E+06 & 87.6\% \\
		\hline
1000	 &0.0384584712 & 7.596407E-07 & 8.7157369E-04 & 11.679 & 8.0108 &	 1.13E+06 & 1.64E+06 & 83.3\% \\
		\hline
1500 &0.0381648481 & 9.266951E-07 & 9.6265025E-04 & 15.490 & 11.79 & 6.97E+05 & 9.15E+05 & 87.7\% \\
		\hline
10000&0.0379344494 & 7.686582E-08 & 2.7724685E-04 & 83.403 & 79.67 & 1.56E+06 & 1.63E+06 & 87.5\% \\
		\hline
100000&0.0381125559&	7.0258327E-09 &8.3820240E-05 & 824.3 & 820.34 & 1.73E+06 & 1.74E+06	 & 87.6\% \\
  		\hline
	\end{tabular}
	\caption{The case $\kappa_s=0.9$, $\kappa_t = 1,$ with importance sampling }
	\label{Tableau_ResultatStat_EP_Ks09_Kt1}
\end{table}

\begin{figure}[!ht]
\centering
	\includegraphics[scale=0.5]{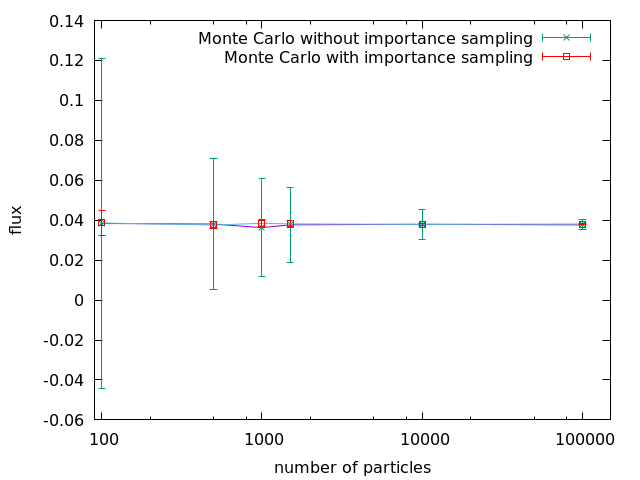}
	\caption{Plot of the results shown in Table~\ref{Tableau_ResultatStat_MC_Ks09_Kt1} and
      Table~\ref{Tableau_ResultatStat_EP_Ks09_Kt1}, that is, $\kappa_s=0.9$ and $\kappa_t=1$.}
	\label{ResultatStationnaireSS09Sa01}
\end{figure}

Table~\ref{Tableau_ResultatStat_MC_Ks09_Kt1} and table~\ref{Tableau_ResultatStat_EP_Ks09_Kt1} show the results for
this test. The last column gives the proportion of particles reaching the inner sphere $r=R_0$. We note that very
few of them reach the inner sphere without the importance sampling method. On the contrary, an important proportion
(almost $90\%$) reach it when the importance sampling is applied. The figure of merit (F.O.M) is computed according to
formula~\eqref{eq:FOM}. When applying the importance sampling method, we provide two execution times, and therefore
two values for the F.O.M. The first one includes the computation of the importance function, which is not
meaningful from a statistical viewpoint, \red{although it is from a computational cost viewpoint}. On the contrary, the second value (Time2, and F.O.M2), do not include it,
and therefore give a clear meaning to the statistical efficiency of the method. As expected, when the number of
particles grows, these two F.O.M are very close to each other. In the case of a small number of particles, the
importance sampling method is less efficient (the F.O.M is increased only by a factor $20$) because the calculation of the importance function is too expensive compared to the Monte
Carlo method. When a large number of particles is used, however, the method is much more efficient, and we see that
the F.O.M is more than $100$ times better than without importance sampling. Figure
\ref{ResultatStationnaireSS09Sa01} shows an important variance reduction, for any number of particles. 

\subsubsection{Case $\kappa_s=0.1$, $\kappa_t=1$}

\begin{table}[!ht]
  \begin{center}
	\begin{tabular}{|c|c|c|c|c|c|c|}
  		\hline
  		N & Flux & Variance & Standard deviation  & Time  & F.O.M & P(N)\\  
  		\hline
      100 & 0.0153756 & 4.3191E-04 & 0.0207825 & 0.0867 & 267018& 0.7\% \\
		\hline
      500 & 0.0248605 & 6.2130E-05 & 0.00788227 & 0.436224 & 368968 & 1.2\%\\
		\hline
      1000 & 0.0220205 & 2.4551E-05 &0.00495486 & 0.847814 & 480437 & 1.04\%\\
		\hline
      1500 & 0.0217672 & 1.0032E-05 &0.00316728 & 1.32612 & 751702 & 1.02\%\\
		\hline
      10000 & 0.0181373 & 4.0237E-06 & 0.00200593 & 8.01809 & 309955 & 0.85\%\\
		\hline
      100000&0.0189385 &1.4290E-07 &3.78026E-04 & 85.582 & 817663 & 0.89\%\\
		\hline
	\end{tabular}    
  \end{center}
	\caption{The case $\kappa_s=0.1$, $\kappa_t = 1,$ without importance sampling.}
	\label{Tableau_ResultatStat_MC_Ks01_Kt1}
\end{table}

\begin{table}[!ht]
	\begin{tabular}{|c|c|c|c|c|c|c|c|c|}
  		\hline
  		N & Flux & Variance & Standard deviation & Time1 & Time2 & F.O.M1 & F.O.M2 & P(N)\\  
  		\hline
      100 & 0.01924 & 1.0554E-06& 0.001027 & 12.2705&	0.58713& 772150 & 1.6137E+07& 88.9\%\\
  		\hline
      500 &	0.01930 & 3.0405E-07& 5.5140E-04& 14.6477&	3.01585& 2.2454E+06& 1.0906E+07& 89.9\%\\
  		\hline
      1000&	0.01906 & 1.2213E-07& 3.4947E-04& 17.5125&	5.9825 & 4.6757E+06& 1.3687E+07& 89.7\%\\
  		\hline
      1500&	0.01900 & 4.9947E-08& 2.2349E-04& 20.2342&	8.79997& 9.8948E+06& 2.2751E+07& 90\%\\
  		\hline
      10000&0.01906 & 1.4452E-08& 1.2022E-04& 71.0115&	59.3438& 9.7439E+06& 1.1659E+07& 89.8\%\\
  		\hline
      100000&0.01910& 3.3606E-09& 5.7970E-05& 596.423&	584.656& 4.9892E+06& 5.0896E+06& 90\%\\
  		\hline
	\end{tabular}
	\caption{The case $\kappa_s=0.1$, $\kappa_t = 1,$ with importance sampling }
	\label{Tableau_ResultatStat_EP_Ks01_Kt1}
\end{table}
\begin{figure}[!ht]
\centering
	\includegraphics[scale=0.5]{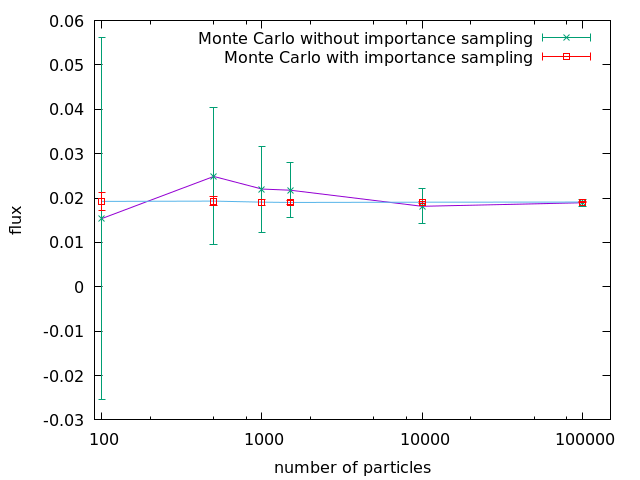}
	\caption{Plot of the results shown in Table~\ref{Tableau_ResultatStat_MC_Ks01_Kt1} and
      Table~\ref{Tableau_ResultatStat_EP_Ks01_Kt1}, that is, $\kappa_s=0.1$ and $\kappa_t=1$.}
	\label{ResultatStationnaireSS01Sa09}
\end{figure}
Table~\ref{Tableau_ResultatStat_MC_Ks01_Kt1} and table~\ref{Tableau_ResultatStat_EP_Ks01_Kt1} show the results for
this test. Here again, the last column gives the proportion of particles reaching the inner sphere $r=R_0$. We note that very
few of them reach the inner sphere without the importance sampling method. On the contrary, an important proportion
(about $90\%$) reach it when the importance sampling is applied. The figure of merit (F.O.M) is computed according to
formula~\eqref{eq:FOM}. 
\red{Here again, we provide two execution times and two F.O.M, the first one including the computation
  of the importance function, the second one excluding it.}
In the case of a small number of particles, the
importance sampling method is less efficient (the F.O.M is increased by a factor $50$) because the calculation of the importance function is too expensive compared to the Monte
Carlo method. When a large number of particles is used, however, the method is much more efficient, and we see that
the F.O.M is more than $150$ times better than without importance sampling. Figure
\ref{ResultatStationnaireSS09Sa01} shows an important variance reduction, for any number of particles. 


\subsection{Unsteady test case}
\label{sec:Cas_Test_Instationnaire}

We use now a test case in which the inner sphere has a radius which depends on time, according to the same law as
in Subsection~\ref{sec:Cas_Test_Validation_Instationnaire}, that is,
\begin{displaymath}
  R_0(t)=0.37625-0.027625t.
\end{displaymath}
As mentioned above, this value of $R_0(t)$ has been derived from three-dimensional simulations of ICF. The incoming
flux imposed on the outer sphere $R_1 = 1$ is constant in time, so we impose a number of particles generated at
each time step. We fix the values of $\kappa_s=0.9$ and $\kappa_t=1$, although other values give the same kind of
results. Here again, we use a uniform mesh in $r$
with $\Delta r= 10^{-2}$, a time step $\Delta t = 10^{-3}$ and a direction discretization $\Delta \mu = 2\times
10^{-3}$. The output value of the code is the flux at the moving boundary $R_0(t)$, integrated from time $t=0$ to
time $t=10$, which is the final time of the simulation. 

The results are presented in Table~\ref{Tableau_Resultat_MC_Ks09_Kt1} et Table
\ref{Tableau_Resultat_MC_EP_Ks09_Kt1}, in which the first column $N$ is the number of particles generated at each
time step. Therefore, the total number of particles generated in the simulation is $N\times T/\Delta t = N\times
10^4.$ As above, the computation time Time1 includes the computation of the importance function, which is done at
each time step, whereas the computation time Time2 does not. The last column in each table represents the
proportion of particles reaching the target.

\begin{table}[ht]
  \begin{center}
	\begin{tabular}{|c|c|c|c|c|c|c|}
      \hline
      N \footnotemark[1] & Flux & Variance & Standard deviation& Time & F.O.M & P(N) \\  
      \hline
      20 & 0.200042 & 8.31746E-06 & 0.002884 & 292 & 412 & 4.44\% \\
      \hline
      100 & 0.200359	 & 6.7959E-07 & 8.2437E-04 & 1437 & 1024 & 4.43\% \\
      \hline
      500 & 0.200569 & 1.9346E-07 & 4.398E-04  & 7722  & 669 & 4.44\%  \\
      \hline
      1000 & 0.2006273 & 1.1675E-07 & 3.417E-04 & 15441 & 554 & 4.44\% \\
      \hline
	\end{tabular}    
  \end{center}
  \caption{The case $\kappa_s=0.9$, $\kappa_t = 1$ without importance sampling.}
  \label{Tableau_Resultat_MC_Ks09_Kt1}
\end{table}
\footnotetext[1]{N is the number of particles generated at each time step.}

\begin{table}[ht]
	\begin{tabular}{|c|c|c|c|c|c|c|c|c|}
  		\hline
  		N & Flux & Variance & Standard deviation  & Time1 & Time2 & F.O.M1 & F.O.M2 & P(N)\\  
  		\hline
  		20 & 0.200611 & 4.07725E-08 &2.019E-04 & 4947 & 639 & 4957 & 38362 & 86\% \\
  		\hline
  		100 & 0.200660& 1.2872E-08 &1.1346E-04 & 7418  & 3092 &10472 & 25125 & 86\% \\
  		\hline
  		500 & 0.2006824 & 2.3614E-09 &4.869E-05 & 20823 & 15923 &20251  & 26481 & 86\%  \\
  		\hline
  		1000 & 0.2006723& 1.1357E-09 &3.37E-05 & 37153 & 31388 &23699&28052 & 86\% \\
  		\hline
	\end{tabular}
	\caption{The case $\kappa_s=0.9$, $\kappa_t = 1$ with importance sampling.}
	\label{Tableau_Resultat_MC_EP_Ks09_Kt1}
\end{table}

\begin{figure}[!ht]
\centering
	\includegraphics[scale=0.4]{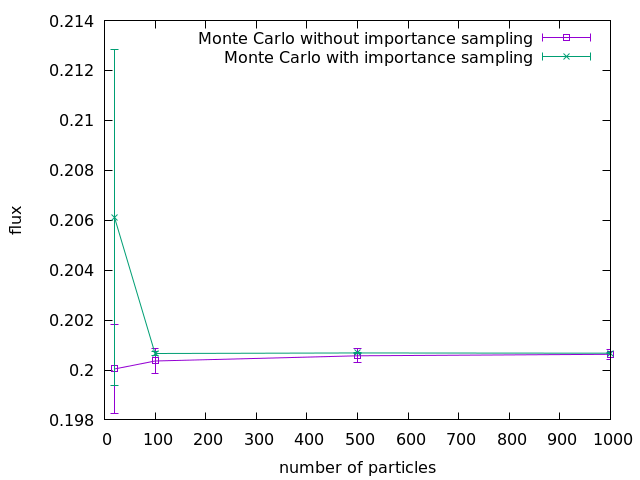}
	\caption{Plot of the results shown in Table~\ref{Tableau_Resultat_MC_Ks09_Kt1} and
      Table~\ref{Tableau_Resultat_MC_EP_Ks09_Kt1}, that is, $\kappa_s=0.9$ et $\kappa_t=1$}
\end{figure}


\section{Conclusion}\label{sec:conclusion}

We have presented in this paper a new method of variance reduction based on importance sampling for the transport
equation in spherical geometry. The importance function is computed as the solution of the adjoint equation, which
is solved numerically. In order to do so, we use an integral equation derived by Siewert and Thomas
\cite{siewert-thomas}, and solve this equation numerically in order to find the first moment (with respect to the
direction $\mu$ of
the importance function). Once this is computed, we apply the method of characteristics to compute the importance
function. Contrary to what has been done in \cite{depinay-these} in a similar context, we do not have an analytical
expression for the importance function. However, it should be noted that the importance function used in
\cite{depinay-these} does not satisfy the correct boundary conditions. Therefore, it is adapted only if boundary
conditions are not of too much importance in the computation at hand. \red{We expect that this is case in ICF
  experiments, in which photons can move a long distance before absorption or scattering. In our method, we are not
  limited by such considerations. }

Numerical tests indicate that the method is efficient, including situations close to the case of inertial
confinement fusion. 

\medskip

\red{In future works, we plan to test this method in situations closer to experiments. In this respect,
  several issues need to be considered:
  \begin{enumerate}
  \item the absorption and scattering coefficients are not constant, contrary to the assumptions we have made
    here. These heterogeneities make the calculation of $I$ (Section~\ref{sec:spherical-case}) much more
    difficult. For instance, the generalization of \eqref{eq:7} might lead to complicated expressions, thereby
    impeding the numerical efficiency of the method. One way (among others) to circumvent this difficulty would
    then be to use a numerical approximation for the computation of $I$ itself. 
  \item The problem is by nature frequency dependent. Although using a grey importance function in such a
    simulation is possible, one should bear in mind that it might prove insufficient. Hence, including a dependence
    of $I$ upon the frequency will probably be an important question to be dealt with. 
  \item The problem we studied here does not take into account interaction with matter. Although we have used
    transient simulations in which this effect is partly represented by the movement of the detector, the
    hydrodynamics of the plasma imply a much richer interaction. This will imply new issues to be considered. 
  \item In relation with the preceding point, a mesh used for an ICF simulation is in general highly
    heterogeneous. This is an additional problem to be considered. 
  \end{enumerate}
}

\appendix

\clearpage
\section{Integral equation on $\phi$}
\label{IntegralEquation}

We give in this Appendix the details of the derivation of the integral equation satisfied by $\phi$. We define
\begin{displaymath}
  \mu_d(r) = -\sqrt{1-\frac{{R_0}^2}{r^2}},
\end{displaymath}
and
\begin{displaymath}
  \phi(r) = r S(r) = \kappa_s r \langle I\rangle.
\end{displaymath}
Equation (\ref{eq:fonction_importance_caracteristiques}) is equivalent to 

\[ \displaystyle
\begin{aligned}
I\left(r,\mu\right) &= \exp\left( \kappa_t\left( r\mu + \sqrt{ r^2\mu^2 - r^2 + {R_0}^2} \right)\right)\mathds{1}_{\left\{ \mu < \mu_d(r)\right\}} + \int_{r\mu}^{R(r,\mu)} S\left(\sqrt{s^2+r^2-r^2\mu^2} \right)\exp\left(\kappa_t \left(r\mu-s\right) \right)ds \\
& = \exp\left( \kappa_t\left( r\mu + \sqrt{ r^2\mu^2 - r^2 + {R_0}^2} \right)\right)\mathds{1}_{\left\{ \mu < \mu_d(r)\right\}} \\
& + \mathds{1}_{\left\{ \mu < \mu_d(r)\right\}}\int_{r\mu}^{-\sqrt{r^2\mu^2-r^2+{R_0}^2}} S\left( \sqrt{s^2 + r^2-r^2\mu^2} \right)\exp\left(\kappa_t \left( r\mu-s\right)\right)  ds \\
& + \mathds{1}_{\left\{ \mu_d(r)<\mu < 0\right\}}\int_{r\mu}^{\sqrt{r^2\mu^2-r^2+{R_1}^2}} S\left( \sqrt{s^2 + r^2-r^2\mu^2} \right)\exp\left(\kappa_t\left( r\mu-s\right)\right) ds \\
& + \mathds{1}_{\left\{ \mu > 0\right\}}\int_{r\mu}^{\sqrt{r^2\mu^2-r^2+{R_1}^2}} S\left( \sqrt{s^2 + r^2-r^2\mu^2} \right)\exp\left(\kappa_t\left( r\mu-s\right)\right)  ds,\\
\end{aligned}
\]
That is,
\[ \displaystyle
\begin{aligned}
I\left(r,\mu\right) &= \exp\left( \kappa_t\left( r\mu + \sqrt{ r^2\mu^2 - r^2 + {R_0}^2} \right)\right)\mathds{1}_{\left\{ \mu < \mu_d(r)\right\}} \\
& + \mathds{1}_{\left\{ \mu < \mu_d(r)\right\}}\int_{r\mu}^{-\sqrt{r^2\mu^2-r^2+{R_0}^2}} S\left( \sqrt{s^2 + r^2-r^2\mu^2} \right)\exp\left(\kappa_t\left( r\mu-s \right)\right)  ds \\
& + \mathds{1}_{\left\{ \mu_d(r)<\mu < 0\right\}}\left( \int_{r\mu}^{0} S\left( \sqrt{s^2 + r^2-r^2\mu^2} \right)\exp\left(\kappa_t\left( r\mu-s\right)\right)ds
\right. \\
&\left.+ \int_{0}^{\sqrt{r^2\mu^2-r^2+{R_1}^2}} S\left( \sqrt{s^2 + r^2-r^2\mu^2} \right)\exp\left(\kappa_t\left( r\mu-s \right)\right) ds\right)    \\
& + \mathds{1}_{\left\{ \mu > 0\right\}}\int_{r\mu}^{\sqrt{r^2\mu^2-r^2+{R_1}^2}} S\left( \sqrt{s^2 + r^2-r^2\mu^2} \right)\exp\left(\kappa_t \left( r\mu-s \right)\right)  ds.\\
\end{aligned}
\]
Hence, changing variables $r' = \sqrt{s^2+r^2-r^2\mu^2}$, we find
\[ \displaystyle
\begin{aligned}
I\left(r,\mu\right) &= \exp\left( \kappa_t\left( r\mu + \sqrt{ r^2\mu^2 - r^2 + {R_0}^2} \right)\right)\mathds{1}_{\left\{ \mu < \mu_d(r)\right\}} \\
& + \mathds{1}_{\left\{ \mu < \mu_d(r)\right\}} \int_{R_0}^{r} r'S\left( r'\right)\exp\left( \kappa_t\left( r\mu+\sqrt{r^2\mu^2-r^2+r'^2} \right) \right) \frac{dr'}{\sqrt{r^2\mu^2-r^2+r'^2}} \\
& + \mathds{1}_{\left\{ \mu_d(r)<\mu < 0\right\}} \left( \int_{r\sqrt{1-\mu^2}}^{r} r'S\left( r'\right)\exp\left( \kappa_t\left( r\mu+\sqrt{r^2\mu^2-r^2+r'^2} \right) \right) \frac{dr'}{\sqrt{r^2\mu^2-r^2+r'^2}}  \right.\\
& + \left.\int_{r\sqrt{1-\mu^2}}^{R_1} r'S\left( r'\right)\exp\left( \kappa_t\left( r\mu-\sqrt{r^2\mu^2-r^2+r'^2} \right) \right) \frac{dr'}{\sqrt{r^2\mu^2-r^2+r'^2}} \right) \\
& + \mathds{1}_{\left\{\mu > 0\right\}}\int_{r}^{R_1} r'S\left( r'\right)\exp\left( \kappa_t\left( r\mu-\sqrt{r^2\mu^2-r^2+r'^2} \right) \right) \frac{dr'}{\sqrt{r^2\mu^2-r^2+r'^2}}. \\
\end{aligned}
\]
Next, we integrate with respect to $\mu$:
\[ \displaystyle
\begin{aligned}
\int_{-1}^{1}I\left(r,\mu\right)d\mu &= \int_{-1}^{1}\exp\left( \kappa_t\left( r\mu + \sqrt{ r^2\mu^2 - r^2 + {R_0}^2} \right)\right)\mathds{1}_{\left\{ \mu < \mu_d(r)\right\}} d\mu \\
& + \int_{-1}^{1}\mathds{1}_{\left\{ \mu < \mu_d(r)\right\}} \int_{R_0}^{r} r'S\left( r'\right)\exp\left( \kappa_t\left( r\mu+\sqrt{r^2\mu^2-r^2+r'^2} \right) \right) \frac{dr'}{\sqrt{r^2\mu^2-r^2+r'^2}} d\mu\\
& + \int_{-1}^{1} \mathds{1}_{\left\{ \mu_d(r)<\mu < 0\right\}} \left( \int_{r\sqrt{1-\mu^2}}^{r} r'S\left( r'\right)\exp\left( \kappa_t\left( r\mu+\sqrt{r^2\mu^2-r^2+r'^2} \right) \right) \frac{dr'}{\sqrt{r^2\mu^2-r^2+r'^2}}  \right.\\
& + \left.\int_{r\sqrt{1-\mu^2}}^{R_1} r'S\left( r'\right)\exp\left( \kappa_t\left( r\mu-\sqrt{r^2\mu^2-r^2+r'^2} \right) \right) \frac{dr'}{\sqrt{r^2\mu^2-r^2+r'^2}} \right) d\mu\\
& + \int_{-1}^{1} \mathds{1}_{\left\{\mu > 0\right\}}\int_{r}^{R_1} r'S\left( r'\right)\exp\left( \kappa_t\left( r\mu-\sqrt{r^2\mu^2-r^2+r'^2} \right) \right) \frac{dr'}{\sqrt{r^2\mu^2-r^2+r'^2}}d\mu. \\
\end{aligned}
\]
We deal with each term seperately.
\begin{itemize}
\item[-]Let \[\displaystyle I_0(r) = \int_{-1}^{\mu_d(r)}\exp\left( \kappa_t\left( r\mu + \sqrt{ r^2\mu^2 - r^2 + {R_0}^2} \right)\right) d\mu.\]
\end{itemize}
Setting
\[ \theta = r\mu + \sqrt{ r^2\mu^2 - r^2 + {R_0}^2}, \] we obtain
\[ \displaystyle
\boxed{I_0(r) = \frac{1}{2r} \left( \left[\frac{1}{\kappa_t}  \exp\left(\kappa_t \theta \right) \right]_{-r+R_0}^{-\sqrt{r^2-{R_0}^2}} + \left( {R_0}^2-r^2 \right) \left[ \kappa_t\text{Ei} \left(\kappa_t\theta \right) - \frac{\exp\left( \kappa_t\theta\right)}{\theta } \right]_{-r+R_0}^{-\sqrt{r^2-{R_0}^2}} \right).}
\]
\begin{itemize}
\item[-]Let
\[ \displaystyle 
\begin{aligned}
I_1(r) &= \int_{-1}^{1}\mathds{1}_{\left\{ \mu < \mu_d(r)\right\}} \int_{R_0}^{r} r'S\left(r'\right)\exp\left( \kappa_t\left( r\mu+\sqrt{r^2\mu^2-r^2+r'^2} \right) \right) \frac{dr'}{\sqrt{r^2\mu^2-r^2+r'^2}} d\mu \\
&=  \int_{R_0}^{r} r'S(r') \left( \int_{-1}^{\mu_d(r)} \exp\left( \kappa_t\left( r\mu + \sqrt{ r^2\mu^2 - r^2 + {r'}^2} \right)\right) \frac{d\mu}{\sqrt{r'^2-r^2 +r^2\mu^2}} \right) dr'
\end{aligned}
\]
\end{itemize}
We change variables as follows:
\[ \nu = r\mu + \sqrt{ r^2\mu^2 - r^2 + {r'}^2}, \] and get
\[ \displaystyle
\begin{aligned}
 I_1(r) & = \int_{R_0}^{r} r'S(r') \left( \int_{-r+r'}^{-\sqrt{r^2-{R_0}^2}+\sqrt{r'^2-{R_0}^2}} \exp\left( \kappa_t \nu\right) \frac{d\nu}{r\nu} \right)  dr' \\
& = \frac{1}{r} \int_{R_0}^{r} r'S(r') \bigg[ \text{Ei}\left( \kappa_t \nu \right) \bigg]_{-r+r'}^{-\sqrt{r^2-{R_0}^2}+\sqrt{r'^2-{R_0}^2}} dr',
\end{aligned}
\]
whence 
\[ \displaystyle
\boxed{I_1(r) = \frac{1}{r} \int_{R_0}^{r} r'S(r') \bigg[ \text{Ei}\left( \kappa_t \left( -\sqrt{r^2-{R_0}^2}+\sqrt{r'^2-{R_0}^2}\right)  \right) - \text{Ei}\left( \kappa_t \left(-r+r' \right)\right) \bigg] dr'.}
\]
\begin{itemize}
\item[-]Let
\[ \displaystyle 
\begin{aligned}
I_2(r) &= \int_{-1}^{1}\mathds{1}_{\left\{ \mu_d(r)<\mu < 0\right\}}  \int_{r\sqrt{1-\mu^2}}^{r} r'S\left( r'\right)\exp\left( \kappa_t\left( r\mu+\sqrt{r^2\mu^2-r^2+r'^2} \right) \right) \frac{dr'}{\sqrt{r^2\mu^2-r^2+r'^2}}d\mu  \\
&=  \int_{R_0}^{r} r'S(r') \left( \int_{\mu_d(r)}^{-\sqrt{1-\frac{r'^2}{r^2}}} \exp\left( \kappa_t\left( r\mu + \sqrt{ r^2\mu^2 - r^2 + {r'}^2} \right)\right) \frac{d\mu}{\sqrt{r'^2-r^2 +r^2\mu^2}} \right) dr'
\end{aligned}
\]
\end{itemize}
Changing variables according to 
\[ \nu = r\mu + \sqrt{ r^2\mu^2 - r^2 + {r'}^2}, \] we find
\[ \displaystyle
\begin{aligned}
 I_2(r) & = \int_{R_0}^{r} r'S(r') \left( \int^{-\sqrt{r^2-r'^2}}_{-\sqrt{r^2-{R_0}^2} +\sqrt{r'^2-{R_0}^2}} \exp\left( \kappa_t \nu\right) \frac{d\nu}{r\nu} \right)  dr' \\
& = \frac{1}{r} \int_{R_0}^{r} r'S(r') \bigg[ \text{Ei}\left( \kappa_t \nu \right) \bigg]^{-\sqrt{r^2-r'^2}}_{-\sqrt{r^2-{R_0}^2}+\sqrt{r'^2-{R_0}^2}} dr',
\end{aligned}
\]
thus
\[ \displaystyle
\boxed{I_2(r) = \frac{1}{r} \int_{R_0}^{r} r'S(r') \bigg[\text{Ei}\left( \kappa_t \left(-\sqrt{r^2-r'^2} \right)\right) -  \text{Ei}\left( \kappa_t \left( -\sqrt{r^2-{R_0}^2}+\sqrt{r'^2-{R_0}^2}\right)  \right) \bigg] dr'.}
\]
\begin{itemize}
\item[-]Let
\[ \displaystyle 
\begin{aligned}
I_3(r) &= \int_{-1}^{1}\mathds{1}_{\left\{ \mu_d(r)<\mu < 0\right\}}  \int_{r\sqrt{1-\mu^2}}^{R_1} r'S\left( r'\right)\exp\left( \kappa_t\left( r\mu-\sqrt{r^2\mu^2-r^2+r'^2} \right) \right) \frac{dr'}{\sqrt{r^2\mu^2-r^2+r'^2}}  \\
&=  \int_{R_0}^{R_1} r'S(r') \left( \int_{\mu_d(r)}^{-\sqrt{1-\frac{r'^2}{r^2}}\mathds{1}_{\left\{ r'<r \right\}}} \exp\left( \kappa_t\left( r\mu - \sqrt{ r^2\mu^2 - r^2 + {r'}^2} \right)\right) \frac{d\mu}{\sqrt{r'^2-r^2 +r^2\mu^2}} \right) dr'
\end{aligned}
\]
\end{itemize} 
Changing variables by setting
\[ \nu = r\mu - \sqrt{ r^2\mu^2 - r^2 + {r'}^2}, \] we infer
\[ \displaystyle
\begin{aligned}
 I_3(r) & = \int_{R_0}^{R_1} r'S(r') \left( -\int^{-\sqrt{|r^2-r'^2|}}_{-\sqrt{r^2-{R_0}^2} - \sqrt{r'^2-{R_0}^2}} \exp\left( \kappa_t \nu\right) \frac{d\nu}{r\nu} \right)  dr' \\
& = \frac{1}{r} \int_{R_0}^{R_1} r'S(r') \bigg[ \text{Ei}\left( \kappa_t \nu \right) \bigg]_{-\sqrt{|r^2-r'^2|}}^{-\sqrt{r^2-{R_0}^2}-\sqrt{r'^2-{R_0}^2}} dr'.
\end{aligned}
\]
Hence,
\[ \displaystyle
\boxed{I_3(r) = \frac{1}{r} \int_{R_0}^{R_1} r'S(r') \bigg[\text{Ei}\left( \kappa_t \left( -\sqrt{r^2-{R_0}^2}-\sqrt{r'^2-{R_0}^2}\right)  \right) - \text{Ei}\left( \kappa_t \left(-\sqrt{|r^2+r'^2|} \right)\right)  \bigg] dr'.}
\]
\begin{itemize}
\item[-]Let
\[ \displaystyle 
\begin{aligned}
I_4(r) &= \int_{-1}^{1} \mathds{1}_{\left\{ \mu > 0\right\}}  \int_{r}^{R_1} r'S(r') \exp\left( \kappa_t\left( r\mu - \sqrt{ r^2\mu^2 - r^2 + {r'}^2} \right)\right) \frac{dr'}{\sqrt{r'^2-r^2 +r^2\mu^2}}d\mu \\
&=  \int_{r}^{R_1} r'S(r') \left( \int_{0}^{1} \exp\left( \kappa_t\left( r\mu - \sqrt{ r^2\mu^2 - r^2 + {r'}^2} \right)\right) \frac{d\mu}{\sqrt{r'^2-r^2 +r^2\mu^2}} \right) dr'
\end{aligned}
\]
\end{itemize}
We change variables as follows:
\[ \nu = r\mu - \sqrt{ r^2\mu^2 - r^2 + {r'}^2}, \] and get
\[ \displaystyle
\begin{aligned}
 I_4(r) & = \int_{r}^{R_1} r'S(r') \left( -\int^{r-r'}_{-\sqrt{r'^2-{r}^2}} \exp\left( \kappa_t \nu\right) \frac{d\nu}{r\nu} \right)  dr' \\
& = \frac{1}{r} \int_{r}^{R_1} r'S(r') \bigg[ \text{Ei}\left( \kappa_t \nu \right) \bigg]_{r-r'}^{-\sqrt{r'^2-{r}^2}} dr',
\end{aligned}
\]
from which we deduce
\[ \displaystyle
\boxed{I_4(r) = \frac{1}{r} \int_{r}^{R_1} r'S(r') \bigg[ \text{Ei}\left( \kappa_t \left( -\sqrt{r'^2-{r}^2} \right)  \right) - \text{Ei}\left( \kappa_t \left(r-r' \right)\right) \bigg] dr'.}
\]
Collecting all the above results, we have
\[ \displaystyle
\begin{aligned}
\int_{-1}^{1} I\left(r,\mu\right) d\mu & = I_0(r) + I_1(r) + I_2(r) + I_3(r) + I_4(r) \\
& = I_0(r) + \frac{1}{r} \int_{R_0}^{R_1} r'S(r') \bigg[ \text{Ei}\left( \kappa_t \left(-\sqrt{r^2-{R_0}^2} -\sqrt{r'^2-{R_0}^2} \right)  \right) - \text{Ei}\left( \kappa_t \left(-|r-r'| \right)\right) \bigg] dr'. \\
\end{aligned}
\]
Recalling the definition of $\phi$, we get
\[ \displaystyle
\begin{aligned}
\phi(r) & = \frac{\kappa_s r}{2} \left( I_0(r) + \frac{1}{r} \int_{R_0}^{R_1} \phi(r') \bigg[ \text{Ei}\left( \kappa_t \left(-\sqrt{r^2-{R_0}^2} -\sqrt{r'^2-{R_0}^2} \right)  \right) - \text{Ei}\left( \kappa_t \left(-| r-r'| \right)\right) \bigg] dr' \right),
\end{aligned}
\]
hence
\[ \displaystyle
\boxed{ 
\begin{aligned}
\phi(r) = \frac{\kappa_s}{4} & \left( \frac{1}{\kappa_t} \bigg[ \exp\left(\kappa_t \theta \right) \bigg]_{-r+R_0}^{-\sqrt{r^2-{R_0}^2}} + \left({R_0}^2 - r^2 \right)\left[ \kappa_t \Ei \left( \kappa_t \theta \right) - \frac{\exp\left( \kappa_t \theta \right)}{\theta} \right]_{-r+R_0}^{-\sqrt{r^2-{R_0}^2}}\right) \\
& + \frac{\kappa_s}{2}\int_{R_0}^{R_1} \phi(r') \bigg[ \Ei \left( \kappa_t \left(-\sqrt{r^2-{R_0}^2} - \sqrt{{r'}^2-{R_0}^2} \right) \right) - \Ei \big( \kappa_t\left(-| r -r' | \right)\big)   \bigg] dr'.
\end{aligned}
}
\]
\clearpage

\bibliographystyle{plain}
\bibliography{MainArticle}

\end{document}